\documentclass[reqno]{amsart}
\usepackage{graphicx}
\usepackage{enumerate}
\usepackage{color}
\newtheorem{thm}{Theorem}[section]
\newtheorem{cor}[thm]{Corollary}
\newtheorem{lem}[thm]{Lemma}
\newtheorem{example}[thm]{Example}
\newtheorem{prop}[thm]{Proposition}
\theoremstyle{definition}
\newtheorem{defn}[thm]{Definition}
\theoremstyle{remark}
\newtheorem{rem}[thm]{Remark}

\numberwithin{equation}{section}

\begin{document}

\title{Lie superalgebras with some homogeneous structures}%

\author{Imen Ayadi , Hedi Benamor , Sa\"{i}d Benayadi}%
\address{Laboratoire de Math\'ematiques et Applications de Metz, CNRS-UMR 7122, Universit\'e Paul Verlaine-Metz,  Ile du Saulcy, F-57045 Metz cedex 1, France.}%
\email{imen.ayadi@univ-metz.fr,benamor@univ-metz.fr,benayadi@univ-metz.fr}%


\begin{abstract}
We generalize to the case of Lie superalgebras the classical symplectic double extension of symplectic Lie algebras introduced in \cite{aubert1996}. We use this concept to give an inductive description of nilpotent homogeneous-symplectic Lie superalgebras. Several examples are included to show the existence of homogeneous quadratic symplectic Lie superalgebras other than even-quadratic even-symplectic considered in \cite{quadsymp}. We study the structures of even (resp. odd)-quadratic  odd (resp. even)-symplectic Lie superalgebras and odd-quadratic odd-symplectic Lie superalgebras and we give its inductive descriptions in terms of quadratic generalized double extensions and odd quadratic generalized double extensions. This study complete the inductive descriptions of homogeneous quadratic  symplectic Lie superalgebras started in \cite{quadsymp}.  Finally, we generalize to the case of homogeneous quadratic symplectic Lie superargebras some relations between even-quadratic even-symplectic Lie superalgebras and Manin superalgebras established in \cite{quadsymp}.
\end{abstract}
\maketitle


\section{Introduction}

In this work, we consider finite dimensional Lie superalgebras over an algebraically closed field ${\mathbb K}$ of characteristic zero. A homogeneous quadratic symplectic Lie superalgebras is a Lie superalgebra endowed with both a homogeneous invariant non-degenerate supersymmetric bilinear form and a homogeneous non-degenerate 2-cocycle of its scalar cohomology. In the non-graded case quadratic symplectic Lie algebras appeared in \cite{aubert1996},\cite{bajo2007}. In particular, in \cite{aubert1996} the quadratic symplectic Lie algebras was studied by using (generalized) symplectic double extension. Next, in \cite{bajo2007}, it was given  new inductive descriptions of quadratic symplectic Lie algebras. The first description was obtained by rewriting the symplectic double extension in terms of quadratic double extension. The second description was obtained by means of the concept of double extension of Manin algebras. Recently in \cite{quadsymp}, a generalization to the case of Lie superalgebra of the quadratic symplectic double extension was done and an inductive description of quadratic symplectic Lie superalgebras  was obtained by using  this notion. Moreover, in the same paper, an other inductive description of quadratic symplectic Lie superalgebras was also obtained by introducing the generalized double extension of Manin superalgebras. The work in \cite{quadsymp} lead us to ask two questions.\\

 1. Does there exist homogeneous quadratic symplectic Lie superalgebras other than even-quadratic even-symplectic Lie superalgebras considered in \cite{quadsymp}?\\

 2. Can we complete the study of the homogeneous quadratic symplectic Lie superalgebras started in \cite{quadsymp}?\\

In this paper, we give an affirmative answer for these questions. Namely, first we prove the existence of homogeneous quadratic symplectic Lie superalgebras by giving some interesting examples of these superalgebras. Next we study the structure of homogeneous-symplectic Lie superalgebras and Lie superalgebras with simultaneously homogeneous-quadratic homogeneous-symplectic structures and we give inductively its descriptions. The main tool used for our purpose is to develop some concept of double extensions. This concept of double extension was introduced by Medina and Revoy in \cite{medina1985} in order to give inductive descriptions of quadratic Lie algebras. Other descriptions of Lie algebras are introduced in \cite{bord1997},\cite{lu2007}. A generalization of the concept of  double extension to the case of Lie superalgebras was introduced in many work \cite{albuquerque2009},\cite{bajo2007generalized},\cite{benamor1999double}, in order to study homogeneous-quadratic Lie superalgebras. In particular, in \cite{albuquerque2009}, it has been proved that any perfect Lie superalgebra with reductive even part can not admits simultaneously even-quadratic and odd-quadratic structures. In the second section of this paper, more precisely in Proposition \ref{prop-abelian-lie-sup-even-quad-odd-quad}, we improve this result by proving that only abelian Lie superalgebras admitting simultaneously even-quadratic and odd-quadratic structures. Moreover, we prove that is not the case for homogeneous-symplectic Lie superalgebras by giving some examples (Example \ref{example-nilpotent} and Example \ref{famille-exemples}) of non-abelian Lie superalgebras which admit simultaneously even-symplectic and odd-symplectic structures. Next, in the third section we study the structure of homogeneous-symplectic Lie superalgebras by introducing the (generalized) double extensions of homogeneous-symplectic Lie superalgebras by a one-dimensional Lie superalgebra and we give an inductive description of nilpotent homogeneous-symplectic Lie superalgebras. The section 4 and 5 are dedicated respectively to give inductive description of odd-quadratic Lie superalgebras with homogeneous-symplectic strucutres and inductive description of even-quadratic odd-symplectic Lie superalgebras. To give these inductive descriptions, we need to define a now concept of generalized double extensions and which are the even-quadratic generalized double extension by the abelian two-dimensional Lie superalgebra and the odd-quadratic generalized double extension by the abelian two-dimensional Lie superalgebra. Finally, in section 6, we generalize to the case of homogeneous quadratic symplectic Lie superargebras some relations between even-quadratic even-symplectic Lie superalgebras and Manin superalgebras established in \cite{quadsymp}. We present some concepts of generalized double extensions of special homogeneous-symplectic homogeneous-Manin superalgebras in order to give other inductive descriptions of homogeneous quadratic symplectic Lie superalgebras.\\

\section{Basic definitions and examples}

\noindent\textbf{1.1 Basic definitions and preliminary results}
\\

\begin{defn} Let ${\frak g}$ be a Lie superalgebra. A bilinear form $\theta$ on ${\frak g}$ is
\begin{enumerate}
   \item odd if $\theta({\frak g}_{\alpha}, {\frak g}_{\alpha}) = \{0\},\ \alpha \in {\mathbb Z}_2$.
   \item even if $\theta({\frak g}_{\alpha}, {\frak g}_{\beta}) = \left\{0\right\}$, where $(\alpha, \beta)\in {\mathbb Z}_2\times {\mathbb Z}_2$ with $\alpha\neq \beta$.
  \item supersymmetric if $\theta(X, Y) = (-1)^{\mid x\mid \mid y\mid} \theta(Y, X),\quad {\forall}_{X \in {\frak g}_{\mid x\mid}, Y \in {\frak g}_{\mid y\mid}}$.
  \item skew-supersymmetric if $\theta(X, Y) = - (-1)^{\mid x\mid \mid y\mid} \theta(Y, X),\quad {\forall}_{X \in {\frak g}_{\mid x\mid},  Y \in {\frak g}_{\mid y\mid}}$.
   \item invariant if $\theta([X, Y], Z) = \theta(X, [Y, Z]),\ {\forall}_{X, Y, Z \in {\frak g}}$.
  \item non-degenerate if $X \in {\frak g}$ satisfies $\Bigl[\theta (X, Y) = 0,\ \forall Y \in {\frak g}\Bigl]$ then $X = 0$.
\end{enumerate}
\end{defn}
\begin{defn}
Let ${\frak g}$ be a Lie superalgebra and $\omega$ a homogeneous bilinear form of degree $\alpha$ on ${\frak g}$ which satisfies the two following conditions:
\begin{enumerate}
\item $\omega$ is a skew-supersymmetric bilinear form,
\item $(-1)^{\mid z\mid \mid x\mid}\omega(X, [Y, Z]) + (-1)^{\mid x\mid \mid y\mid}\omega(Y, [Z, X]) + (-1)^{\mid y\mid \mid z\mid}\omega(Z, [X, Y]) = 0$, ${\forall}_{X\in {\frak g}_{\mid x\mid},Y\in {\frak g}_{\mid y\mid},Z\in {\frak g}_{\mid z\mid}}$.
\end{enumerate}
We say that $\omega$ is a homogeneous scalar $2-$cocycle of ${\frak g}$ of degree $\alpha$.
\end{defn}
 \begin{defn}
 A Lie superalgebra ${\frak g}$ is:
 \begin{enumerate}
 \item \textsl{even (resp. odd)-symplectic} if there exists an even (resp. odd) bilinear form $\omega$ on ${\frak g}$ such that $\omega$ is a non-degenerate scalar $2-$cocycle. It is denoted by $({\frak g}, \omega)$ and $\omega$ is called an even (resp. odd)-symplectic form on ${\frak g}$.
 \item \textsl{even (resp.odd)-quadratic} if there exists an even (resp. odd) bilinear form $B$ on ${\frak g}$ such that $B$ is a supersymmetric, non-degenerate and invariant. It is denoted by $({\frak g}, B)$ and $B$ is called an even (resp. odd)-invariant scalar product on ${\frak g}$.\\
 \end{enumerate}
 \end{defn}

The following definition giving in \cite{albuquerque2009} in order to study odd-quadratic Lie superalgebras, will be very useful in this paper.\\

\begin{defn}
Let ${\frak g}$ be a Lie superalgebra. We denote by $P({\frak g})={V}_{\bar{0}}\oplus {V}_{\bar{1}}$ the ${\mathbb Z}_2$-graded vector space obtained from ${\frak g}$ with graduation defined by
$${V}_{\bar{0}}:= {\frak g}_{\bar{1}}\ \ \mbox{and}\ \ \ {V}_{\bar{1}}:= {\frak g}_{\bar{0}}.$$
Moreover, if we consider ${\frak g}^*$ the dual space of ${\frak g}$, then $P({\frak g}^*)$ and ${\frak g}^*$ are equal as ${\mathbb Z}_2$-graded vector space such that
$${V}^*_{\bar{0}}= {\frak g}^*_{\bar{1}}\ \ \ \mbox{and}\ \ \ {V}^*_{\bar{1}}= {\frak g}^*_{\bar{0}}.$$\\
\end{defn}

When we study Lie superalgebra with homogeneous structures, a natural question arises: what about the compatibility of these structures? An answer of this question, in the case of homogeneous-quadratic perfect Lie superalgebra with reductive even part, is given in \cite{albuquerque2009} where odd-quadratic Lie superalgebras have been studied. In the following proposition, we improve the result obtained in \cite{albuquerque2009}.\\


\begin{prop}\label{prop-abelian-lie-sup-even-quad-odd-quad}
Let $({\frak g}, {\left[,\right]}_{\frak g})$ be a Lie superalgebra. ${\frak g}$ admits simultaneously an even-quadratic structure and an odd-quadratic structure if and only if ${\frak g}$ is an abelian Lie superalgebra such that $dim {\frak g}_{\bar{0}}= dim {\frak g}_{\bar{1}}$.\\
\end{prop}

\noindent\underline{\sl Proof} : Suppose that ${\frak g}$ is provided with an even-quadratic structure and odd-quadratic-structure which are denoted respectively by $B$ and $\bar{B}$. Using the non-degeneracy and the invariance of $B$ and $\bar{B}$, we deduce that the two maps
\vskip 2pt \noindent $\phi: {\frak g}_{\bar{0}} \longrightarrow {\frak g}_{\bar{0}}^*$ and $\psi: {\frak g}_{\bar{1}} \longrightarrow {\frak g}_{\bar{0}}^*$ defined respectively by:
$$\phi(X) =B(X,.), \ \ \mbox{and}\ \ \psi(Y)= \bar{B}(Y,.),\ {\forall}_{X\in {\frak g}_{\bar{0}}, Y\in {\frak g}_{\bar{1}}},$$
are isomorphisms of ${\frak g}_{\bar{0}}$-modules. Consequently, we obtain that ${\frak g}_{\bar{0}}$ and ${\frak g}_{\bar{1}}$ are isomorphic as ${\frak g}_{\bar{0}}$-modules by means of $\Phi:= {\psi}^-\circ \phi$. On the other hand as $({\frak g}_{\bar{1}}, {B}_{\bar{1}}:={B\vert}_{{\frak g}_{\bar{1}}\times {\frak g}_{\bar{1}}})$ is a symplectic vector space such that ${B}_{\bar{1}}$ is ${\frak g}_{\bar{0}}$-invariant then, we deduce that the following bilinear form $$T: {\frak g}_{\bar{0}}\times {\frak g}_{\bar{0}} \longrightarrow {\mathbb K}\ \ \mbox{ defined by}\ \ \ T(X,Y):= B_{\bar{1}}(\Phi(X),\Phi(Y)),$$ is a skew-symmetric, invariant and non-degenerate. So, for $X,Y,Z\in {\frak g}_{\bar{0}}$ we have
\begin{eqnarray*}
T(\left[X,Y\right],Z)&=& T(X,\left[Y,Z\right])= - T(X,\left[Z,Y\right]) = - T(\left[X,Z\right],Y)\\
&=& T(\left[Z,X\right],Y) = T(Z,\left[X,Y\right]) = - T(\left[X,Y\right],Z).
\end{eqnarray*}  \\
Meaning by the non-degeneracy of $B_{\bar{1}}$, that ${\frak g}_{\bar{0}}$ is an abelian Lie algebra. Besides, using the invariance of $\bar{B}$ it comes $\bar{B}(\left[{\frak g}_{\bar{0}}, {\frak g}_{\bar{1}}\right], {\frak g}_{\bar{0}})= \bar{B}({\frak g}_{\bar{1}}, \left[ {\frak g}_{\bar{0}}, {\frak g}_{\bar{0}}\right])= \left\{0\right\}$ and which implies by the non-degeneracy of $\bar{B}$ that $\left[{\frak g}_{\bar{0}}, {\frak g}_{\bar{1}}\right]= \left\{0\right\}$. Finally, by the invariance of $B$, it comes $B(\left[{\frak g}_{\bar{1}}, {\frak g}_{\bar{1}}\right], {\frak g}_{\bar{0}})= B({\frak g}_{\bar{1}}, \left[{\frak g}_{\bar{1}}, {\frak g}_{\bar{0}}\right])= \left\{0\right\}$ and we deduce by the non-degeneracy of $B$ that $\left[{\frak g}_{\bar{1}}, {\frak g}_{\bar{1}}\right]= \left\{0\right\}$. So, we  conclude that ${\frak g}$ is an abelian Lie superalgebra. The fact that $\bar{B}$ is non-degenerate and odd implies that $dim {\frak g}_{\bar{0}}= dim {\frak g}_{\bar{1}}$. $\square$
\\

The following proposition generalize to homogeneous quadratic case and homogeneous symplectic case the
characterization established in \cite{quadsymp}.

\begin{prop}\label{caracterisation}
Let $({\frak g},B)$ be a homogeneous-quadratic Lie superalgebra.
\begin{enumerate}[(i)]
\item If $({\frak g},B)$ is an even-quadratic Lie super-algebra, then ${\frak g}$ admits an odd-symplectic structure $\omega$ if and only if there is an odd invertible skew-symmetric superderivation $\Delta$ of $({\frak g}, B)$ such that
\begin{equation}
\label{equation-prop-caracterisation} \quad \omega(X, Y) = B\bigl(\Delta(X), Y\bigl), \ \forall X,Y\in {\frak g}.
\end{equation}
\item If $({\frak g},B)$ is an odd-quadratic Lie superalgebra, then ${\frak g}$ admits an odd-symplectic structure $\omega$ if and only if there is an even invertible skew-symmetric superderivation $\Delta$ of $({\frak g}, B)$ such that the equation (\ref{equation-prop-caracterisation}) is hold.
\item If $({\frak g},B)$ is an odd-quadratic Lie superalgebra, then ${\frak g}$ admits an even-symplectic structure $\omega$ if and only if there is an odd invertible skew-symmetric superderivation $\Delta$ of $({\frak g}, B)$ such that the equation (\ref{equation-prop-caracterisation}) is hold.
\end{enumerate}
\end{prop}
\noindent\underline{\sl Proof} : The proof of this proposition is analogous to the Proposition 4.1 of \cite{quadsymp}.   $\square$ \\

\begin{rem}\label{homogeneous-symplectic-case}
Following the previous proposition, we deduce that every even (resp. odd)-quadratic Lie superalgebra that admits odd (resp. even)-symplectic structure it admits also even (resp. odd)-symplectic structure. Indeed, let $({\frak g},B,\omega)$ be an even (resp. odd)-quadratic odd (resp. even)-symplectic Lie superalgebra. By the assertion $(ii)$ of the above proposition, there exist an odd invertible skew-supersymmetric superderivation $\Delta$ of $({\frak g},B)$ such that $\omega(X,Y)= B(\Delta(X),Y)$, ${\forall}_{X,Y\in {\frak g}}$. Consider $\widetilde{\Delta}:= \left[\Delta,\Delta\right]$. Clearly $\widetilde{\Delta}$ is an even invertible skew-supersymmetric superderivation of $({\frak g},B)$. So, by the assertion $(i)$ of the above proposition, we conclude that ${\frak g}$ admits even (resp. odd)-symplectic structure.\\
\end{rem}


Symplectic Lie algebras and quadratic symplectic Lie algebras have been studied in \cite{aubert1996} and \cite{bajo2007} where many examples of these algebras are included. Since our goal is to study homogeneous-symplectic Lie superalgebras and homogeneous quadratic symplectic Lie superalgebras, in this subsection we give some interesting examples of these superalgebras. Among The examples we will give, there are non-abelian Lie superalgebras which admit simultaneously even-symplectic and odd-symplectic structures. This phenomenon is different from homogeneous-quadratic Lie superalgebras. Indeed, it was proved in Proposition \ref{prop-abelian-lie-sup-even-quad-odd-quad} that only abelian Lie superalgebras admit simultaneously even-quadratic and odd-quadratic structures.\\

\begin{example}
Consider the two-dimensional Lie superalgebra ${\frak m}= {\mathbb K}e_0\oplus{\mathbb K}e_1$ such that $\left[e_0,e_1\right]= e_1$. Clearly, that the skew-supersymmetric bilinear form $\omega$ on ${\frak m}$ defined by $\omega(e_0,e_1)= 1$ and $\omega$ is null elsewhere is an odd-symplectic form on ${\frak m}$. Moreover, it is easy to show that every two-dimensional odd-symplectic Lie superalgebra is either the abelian Lie superalgebra or isomorphic to ${\frak m}$.\\

\end{example}
\begin{example}\label{example-otimes}
Let $(A,., B)$ be an odd-symmetric associative super-commutative superalgebra. The $2-$dimensional odd-symmetric associative superalgebras defined in \cite{imen-asso} are examples of this type of superalgebras. In addition, let us  consider $({\frak h},{\left[,\right]}_{\frak h} \omega)$ a symplectic Lie algebra. The ${\mathbb Z}_2$-gradeed vector space ${\frak g}:= {\frak h}\otimes A$ (where ${\frak g}_{\bar{0}}:= {\frak h}\otimes A_{\bar{0}}$ and ${\frak g}_{\bar{1}}:= {\frak h}\otimes A_{\bar{1}}$) provided with the following bracket:
\begin{eqnarray}
\label{bracket-otimes}\left[X\otimes a, Y\otimes b\right]:= {\left[X,Y\right]}_{\frak h}\otimes a.b,
\end{eqnarray}
and the following bilinear form: $$\Omega(X\otimes a, Y\otimes b):= \omega(X,Y)B(a,b), \ \forall X,Y\in {\frak h} \ \mbox{and}\ a,b\in A$$
is an odd-symplectic Lie superalgebra.\\
\end{example}

To give other examples of homogeneous quadratic symplectic Lie superalgebras, we need some notions. We recall that for any Lie superalgebra ${\frak h}$, we can consider:\\

\textbf{(1)} ${\frak g}:= {\frak h}\oplus {\frak h}^*$ the trivial even double extension of ${\frak h}$ (i.e the double extension of $\left\{0\right\}$ by ${\frak h}$). The bracket and the invariant product scalar of ${\frak g}$ are given respectively  by
$${\left[X + f,Y + h\right]}_{\frak g}:= {\left[X,Y\right]}_{\frak h} +  X.h - {(-1)}^{\mid x\mid \mid y\mid}Y.f,$$
$$\mbox{and}\ \ \ B(X+f,Y+h)= f(Y) + {(-1)}^{\mid x\mid \mid y\mid}h(X), \ \forall (X+f)\in {\frak g}_{\mid x\mid}, (Y+h)\in {\frak g}_{\mid y\mid}.\\$$

\textbf{(2)} ${\frak g}:= {\frak h}\oplus P({\frak h}^*)$ the trivial odd double extension of ${\frak h}$ (i.e the odd double extension of $\left\{0\right\}$ by ${\frak h}$). Its bracket is equal to the bracket of the trivial even double extension of ${\frak h}$ and its odd-invariant product scalar is given by
$$B(X+f,Y+h)= f(Y) + h(X), \ \forall (X+f)\in {\frak g}_{\mid x\mid}, (Y+h)\in {\frak g}_{\mid y\mid}.$$
\linebreak
The following lemma, given in \cite{benayadi2003socle}, will be very useful to construct some interesting examples.\\

\begin{lem}\label{lemme-du-socle}
Let ${\frak h}$ be a Lie superalgebra and $({\frak g},B)$ the trivial even or odd double extension of ${\frak h}$. Let $D$ be a homogeneous invertible superderivation of ${\frak h}$ of degree $d$ and denote by $D^*$ the homogeneous endomorphism of ${\frak h}^*$ of degree $d$ defined by $D^*(f)(X) := -{(-1)}^{\alpha d}f \circ D(X)$, $\forall f\in {({\frak h}^*)}_{\alpha}$ and $X\in {\frak h}$. Then, the map $\widetilde{D}$ defined by
$$\widetilde{D}(X+f):= D(X)+D^*(f),\ \forall (X+f)\in {\frak g}_{\mid x\mid},$$
is a homogeneous invertible skew-supersymmetric superderivation of ${\frak g}$ of degree $d$.\\
\end{lem}

Applying remark \ref{homogeneous-symplectic-case}, we deduce that Lie superalgebras obtained in Example \ref{example-nilpotent} and Lie superalgebras obtained in Example \ref{famille-exemples} below admit simultaneously even-symplectic and odd-symplectic structures.\\

\begin{example}\label{example-nilpotent}
Let ${\frak L}:= {\frak L}_{\bar{0}}\oplus {\frak L}_{\bar{1}}$, where ${\frak L}_{\bar{0}}:={\mathbb K}l_0\oplus {\mathbb K}k_0$ and ${\frak L}_{\bar{1}}:= {\mathbb K}l_1\oplus {\mathbb K}k_1$, be a Lie superalgebra with bracket defined by $\left[l_0,l_1\right]=k_1$, $\left[l_1,l_1\right]=k_0$ and null elsewhere. On ${\frak L}$, we define the two linear maps $D$ and $\Delta$ respectively by:
\begin{eqnarray}
D(l_i)&=& l_i,\ \ D(k_i)= 2k_i, \ \mbox{where} \ i\in \left\{0,1\right\},\\
\Delta(l_0)&=& l_1,\ \Delta(l_1)= l_0, \ \Delta(k_0)= 2k_1,\ \Delta(k_1)= k_0.
\end{eqnarray}
By definition, it is clear that $D$ is even and $\Delta$ is odd. Moreover, by a simple computation, we can see easily that $D$ and $\Delta$ are two invertible superderivations on ${\frak L}$.  Consequently, by the Lemma \ref{lemme-du-socle}, we obtain that:
\begin{enumerate}
\item[(i)] $({\frak L}\oplus P({\frak L}^*), B, \widetilde{D})$ is an odd-quadratic odd-symplectic Lie superalgebra.
\item[(ii)]$({\frak L}\oplus P({\frak L}^*), B, \widetilde{\Delta})$ is an odd-quadratic even-symplectic Lie superalgebra.
\item[(iii)] $({\frak L}\oplus {\frak L}^*, B, \widetilde{\Delta})$ is an even-quadratic odd-symplectic Lie superalgebra.\\
\end{enumerate}
\end{example}

\begin{example}\label{famille-exemples}
Let us consider the super-vector space $A:= <e_1,\cdots,e_n, f_1,\cdots, f_n>$, where $n\geq 1$, such that $A_{\bar{0}}= <e_1,\cdots, e_n>$ and $A_{\bar{1}}= <f_1,\cdots, f_n>$. Define the bilinear map $".": A\times A \longrightarrow A$ by
\begin{eqnarray*}
e_i.e_j &:=&\left\lbrace
           \begin{array}{l} e_{i+j}, \ i+j\leq n\\
           0 \ \mbox{elsewhere}
           \end{array}
           \right.\\
e_i. f_j &:=&\left\lbrace
           \begin{array}{l} f_{i+j}, \ i+j\leq n\\
           0 \ \mbox{elsewhere}
           \end{array}
           \right.\\
f_i.f_j &:=& 0, \ \forall i,j\in \left\{1,\cdots,n\right\}.
\end{eqnarray*}
provided with $"."$, it is easy to show that $A$ is an associative super-commutative superalgebra. Moreover, let ${\frak h}$ be a Lie algebra. On the Lie superalgebra ${\frak g}:= {\frak h}\otimes A$ which is defined as Example \ref{example-otimes} and provided with bracket (\ref{bracket-otimes}), we define the following endomorphism $D$ for homogeneous elements by
\begin{eqnarray*}
D(X\otimes e_i)&=& i X\otimes f_i,\ \ \ D(X\otimes f_i)= X\otimes e_i, \ \forall X\in {\frak h}.
\end{eqnarray*}
By a simple computation, we can see that $D$ is an odd invertible superderivation on ${\frak g}$. Consequently, by the Lemma \ref{lemme-du-socle}, we obtain that
\begin{enumerate}
\item[(i)] $({\frak g}\oplus {\frak g}^*,B,\widetilde{D})$ is an even-quadratic odd-symplectic Lie superalgebra,
\item[(ii)] $({\frak g}\oplus P({\frak g}^*),B,\widetilde{D})$ is an odd-quadratic even-symplectic Lie superalgera.\\
\end{enumerate}
\end{example}

\begin{rem}
 Since all four-dimensional Lie superalgebras of Example \ref{famille-exemples} are abelian, we deduce that Lie superalgebras of Example \ref{example-nilpotent} can not be obtained as Lie superalgebras of Example \ref{famille-exemples}.\\
\end{rem}

\section{Homogeneous-symplectic nilpotent Lie superalgebras}

In this section we give inductive descriptions of nilpotent homogeneous-symplectic Lie superalgebras. Our main tool to give these inductive descriptions is to generalize to the case of Lie superalgebras the classical symplectic double extension of symplectic Lie algebras introduced in \cite{aubert1996}. We start with the following lemma.

\begin{lem}\label{orthogonal-of-I-is-ideal}
Let $({\frak g},\omega)$ be a homogeneous-symplectic Lie superalgebra, $I$ a graded ideal of ${\frak g}$ contained in ${\frak z}({\frak g})$ and $J$ its orthogonal with respect to $\omega$. Then, $J$ is a graded ideal of ${\frak g}$.
\end{lem}

\noindent\underline{\sl Proof} : Since $\omega$ is homogeneous (even or odd) and $I$ is graded it comes that $J$ is graded. Using the fact that $I$ is contained in ${\frak z}({\frak g})$ and the fact that $\omega$ is a scalar 2-cocycle of ${\frak g}$, we deduce that $J$ is a graded ideal of ${\frak g}$. \ \ \ \ \ $\square$ \\

\noindent\textbf{2.1 Inductive description of odd-symplectic nilpotent Lie superalgebras}
\\

 We introduce the concept of double extension and the concept of the generalized double extension of an odd-symplectic Lie superalgebra by a one-dimensional Lie superalgebra ${\mathbb K}e$. In the first step, we  obtain the central extension of $P({\mathbb K}e^*)$ by an odd-symplectic Lie superalgebra (we consider $\left\{e^*\right\}$ the dual basis of $\left\{e\right\}$).\\


\begin{prop}\label{extension-centrale}
Let $({\frak g},\omega)$ be an odd-symplectic Lie superalgebra and $D$ a homogeneous superderivation on ${\frak g}$ of degree $\mid e\mid$. Define the bilinear form $\gamma: {\frak g}\times {\frak g}\longrightarrow {\mathbb K}$ for homogeneous elements by
\begin{equation*}
\label{deux-cocycle}\gamma(X, Y) = - \Bigl[\omega\bigl(D(X), Y\bigl) + (-1)^{\mid e\mid \mid x\mid} \omega\bigl(X, D(Y)\bigl)\Bigl].
\end{equation*}
Then the ${\mathbb Z}_2$-graded vector space $P({\mathbb K}e^*) \oplus {\frak g}$ endowed with the following bracket
\begin{equation}
 \quad [a e^*+X, b e^*+Y] :=  \gamma(X, Y) e^* + [X, Y]_{\frak g}, {\forall}_{(a e^*+ X),(b e^* +Y)\in P({\mathbb K}e^*)\oplus{\frak g}}
\end{equation}is a Lie superalgebra. $P({\mathbb K}e^*)\oplus {\frak g}$ is called the central extension of ${\frak g}$ by $P({\mathbb K}e^*)$ by means of $\gamma$.\\
\end{prop}

\noindent\underline{\sl Proof} : By a simple computation where we use that fact that $\omega$ is a scalar 2-cocycle of ${\frak g}$ and the fact that $D$ is a homogeneous superderivation of ${\frak g}$ of degree $\mid e\mid$, we deduce easily that $\gamma$ is a scalar 2-cocycle of ${\frak g}$ and so the result. \ \ \ \ $\square$\\

\begin{thm}\label{g-d-e-odd-symp}
Let $({\frak g},\omega)$ be an odd-symplectic Lie superalgebra, $D$ a homogeneous superderivation on ${\frak g}$ of degree $\mid e\mid$ and $x_0\in {\frak g}_{\bar{0}}$ such that
\begin{eqnarray}
\label{condition-g-d-e-odd-symp}(1 - {(-1)}^{\mid e\mid})(D^2- \frac{1}{2}\left[x_0,.\right])=0\ \ \ \mbox{and}\ \ \ (1 - {(-1)}^{\mid e\mid})D(x_0)=0.
\end{eqnarray}
 Consider the following scalar 2-cocycle of ${\frak g}$ defined by
$$\theta(X,Y)= \omega(\Bigl[D^*({(-1)}^{\mid e\mid \mid x\mid}D+ D^*) + (D + {(-1)}^{\mid e\mid (\mid x\mid + 1)} D^*)D\Bigl](X),Y), \ \ {\forall}_{X,Y\in {\frak g}},$$
where $D^*$ is the adjoint of $D$ with respect to $\omega$. If $\theta$ is a 2-coboundary of ${\frak g}$ (i.e there exist $b_0\in {\frak g}_{\bar{0}}$ such that $\theta(X,Y) = \omega(b_0, [X,Y]),\ {\forall}_{X,Y\in {\frak g}}$), which is the case if $D$ is an odd superderivaion for $b_0= \frac{1}{2}x_0$ then the linear map $\widetilde{D}$ on $P({\mathbb K}e^*)\oplus{\frak g}$ defined by
$$\widetilde{D}(a e^*+X):= D(X) - \omega(b_0, X)e^*,\ {\forall}_{a e^*+X\in P({\mathbb K}e^*)\oplus{\frak g}}$$
is a homogeneous superderivation on $P({\mathbb K}e^*)\oplus{\frak g}$ of degree $\mid e\mid$. Consequently, the ${\mathbb Z}_2$-graded vector space ${\frak t}= P({\mathbb K}e^*) \oplus {\frak g}\oplus {\mathbb K}e$ endowed with the following bracket defined by
\begin{eqnarray*}\label{produitimpaire}
[e, e] &=& \frac{1}{2}(1 - {(-1)}^{\mid e\mid})x_0 \ ; \quad [e, X] = \widetilde{D}(X);\\
 \quad [X, Y] &=& {\left[X,Y\right]}_{\frak g} +\gamma(X, Y) e^* ,\ {\forall}_{X, Y \in {\frak g}}.
\end{eqnarray*}
is a Lie superalgebra. Moreover the  skew-supersymmetric bilinear form $\Omega$ on ${\frak t}$ defined by
\begin{equation}
\Omega_{\vert_{\frak g\times \frak g}} = \omega,\ \ \Omega(e^*, e) = 1, \ \ \ \mbox{and}\ \ \Omega \ \mbox{ is null elsewhere},
\end{equation}
is an odd-symplectic form on ${\frak t}$. $({\frak t}, \Omega)$ is an odd-symplectic Lie superalgebra which we called the double extension (resp. generalized double extension )of $({\frak g}, \omega)$ by the one-dimensional Lie algebra (resp. the one-dimensional Lie superalgebra with null even part) by means of $(D,b_0)$ (resp. $(D,x_0)$).\\
\end{thm}

\noindent\underline{\sl Proof} :  we start by proving that $\widetilde{D}\in {(Der(P({\mathbb K}e^*)\oplus {\frak g}))}_{\mid e\mid}$. Meaning we have to show that
$$\widetilde{D}({\left[X,Y\right]}_{P({\mathbb K}e^*)\oplus{\frak g}})= {\left[\widetilde{D}(X),Y\right]}_{P({\mathbb K}e^*)\oplus{\frak g}} +  {(-1)}^{\mid e\mid \mid x\mid}{\left[X,\widetilde{D}\right]}_{P({\mathbb K}e^*)\oplus{\frak g}}, {\forall}_{X,Y\in {\frak g}}.$$
Analyzing each term, we get:
\begin{eqnarray*}
\widetilde{D}({\left[X,Y\right]}_{P({\mathbb K}e^*)\oplus{\frak g}}) &=& D(\left[X,Y\right]) - \omega(b_0, \left[X,Y\right])e^*,\\
{\left[\widetilde{D}(X),Y\right]}_{P({\mathbb K}e^*)\oplus{\frak g}} +  {(-1)}^{\mid e\mid \mid x\mid}{\left[X,\widetilde{D}\right]}_{P({\mathbb K}e^*)\oplus{\frak g}}&=& \left[D(X),Y\right] +  {(-1)}^{\mid e\mid \mid x\mid} \left[X,D(Y)\right]\\
 &-& \theta(X,Y).
\end{eqnarray*}
Since $D$ is a superderivation of degree $\mid e\mid$ and $\theta(X,Y):= \omega(b_0, \left[X,Y\right])$, we have $\widetilde{D}$ is a homogeneous superderivation of $P({\mathbb K}e^*)\oplus{\frak g}$ of degree $\mid e\mid$. So we conclude that we can consider ${\frak t}= P({\mathbb K}e^*) \oplus {\frak g}\oplus {\mathbb K}e$ the semi-direct product (resp. the generalized semi-direct product) of $P({\mathbb K}e^*) \oplus {\frak g}$ by the one-dimensional Lie algebra (resp. the one-dimensional Lie superalgebra with non even part) ${\mathbb K}e$ by means of $\widetilde{D}$ (resp. $(\widetilde{D},x_0)$) and so the result. \ \ \ \ \ $\square$ \\

\begin{prop}\label{prop-inverse-impaire} Let $({\frak g}, \omega)$ be an odd-symplectic Lie superalgebra such that ${\frak z}({\frak g})\not = \{0\}$, then $({\frak g}, \omega)$ is either a double extension or a generalized double extension of an odd-symplectic Lie superalgebra $({\frak h}, {\omega}_{\frak h})$ by a one-dimensional Lie superalgebra.\\
 \end{prop}

 \noindent\underline{\sl Proof} : Let $e^*$ be a homogeneous non zero element of ${\frak z}({\frak g})$ and $I:={\mathbb K}e^*$ the ideal of ${\frak g}$ generated by $e^*$. By Lemma \ref{orthogonal-of-I-is-ideal}, it comes that $J={I}^{\bot}$ is a graded ideal of ${\frak g}$. Since $\omega$ is odd, it is clear that $I\subseteq J$. As $\omega$ is odd and non-degenerate then there exists $e\in {\frak g}_{\mid e^*\mid + 1}$ such that $\omega(e^*,e)\neq 0$. From the non-degeneracy of $\omega$, there exists a graded vector space ${\frak h}$ of ${\frak g}$ such that $\omega(e^*, {\frak h}) = \{0\} = \omega(e, {\frak h}), J = I \oplus {\frak h}$ and ${\omega}_{\frak h} = \omega_{|_{\frak h\times\frak h}}$ is a non-degenerate form on ${\frak h}$. Since ${\frak h}$ is a graded vector space of ${\frak g}$ contained in the graded ideal $J$ then the bracket of ${\frak g}$ takes the form :
\begin{eqnarray*}
[X, Y] &=& \gamma(X, Y) e^* + \alpha(X,Y),
\end{eqnarray*}
where $\alpha(X,Y)\in {\frak h}$ and $\gamma(X,Y)\in {\mathbb K}$. Since  ${\frak g}$ is a Lie superalgebra, it comes directly that $({\frak h},\alpha)$ is a Lie superalgebra. Besides it not difficult to see that ${\omega}_{\frak h}$ is an odd-symplectic from on ${\frak h}$. Now concerning the rest of bracket in ${\frak g}$ takes the form
$$[e, X] = D(X) + \beta(X)e^*, \ \ \ \mbox{and}\ \ \ [e, e] = \frac{1}{2}( 1 - {(-1)}^{\mid e\mid})(k e^* + x_0), {\forall}_{X\in {\frak h}},$$
where $D(X)\in {\frak h}$, $x_0\in {\frak h}_{0}$ and $\beta(X), k\in {\mathbb K}$. The fact that $\omega$ is a scalar 2-cocycle of ${\frak g}$ yields that $k=0$, $(1 - {(-1)}^{\mid e\mid})(\beta(X) + \frac{1}{2}\omega(x_0,X))=0$ and $$\gamma (X, Y) = - \Bigl[{\omega}_{\frak h}\bigl(D(X), Y\bigl) + (-1)^{\mid e\mid \mid x\mid} {\omega}_{\frak h}\bigl(X, D(Y)\bigl)\Bigl], {\forall}_{X,Y\in {\frak h}}.$$ On the other hand, by the graded identity of Jacobi, it follows that $(1 - {(-1)}^{\mid e\mid}) D(x_0)=0$, $(1 - {(-1)}^{\mid e\mid}) (D^2 - \frac{1}{2}\left[x_0,.\right])=0$, $D$ is a homogeneous superderivation of ${\frak h}$ of degree $\mid e\mid $ and for homogeneous elements $X$ and $Y$ in ${\frak h}$, we have
\begin{equation}
\label{condition-cobord}\beta(\left[X,Y\right])= - {\omega}_{\frak h}(\Bigl[D^*({(-1)}^{\mid e\mid \mid x\mid}D+ D^*) + (D + {(-1)}^{\mid e\mid (\mid x\mid + 1)} D^*)D\Bigl](X),Y).
\end{equation}
 Meaning that if $D$ is an odd superderivation, then the condition (\ref{condition-cobord}) is trivially satisfied whereas if $D$ is an even superderivation, then the condition (\ref{condition-cobord}) implies that the bilinear form
$$\theta(X,Y):= -{\omega}_{\frak h}(\Bigl[ D^*(D+ D^*) + (D + D^*)D\Bigl](X),Y), {\forall}_{X,Y \in {\frak h}},$$
is a scalar 2-coboundary of ${\frak h}$.  Since ${\omega}_{\frak h}$ and $\beta\in {\frak h}^*$, then there exists a unique $z_0\in {\frak h}$ such that $\beta= \omega(z_0,.)$. Now applying Theorem \ref{g-d-e-odd-symp}, we obtain that ${\frak g}$ is either the double extension of $({\frak h},{\omega}_{\frak h})$ by the one-dimensional Lie algebra ${\mathbb K}e$ by means of $(D,z_0)$ or the generalized double extension of $({\frak h},{\omega}_{\frak h})$ by the one-dimensional Lie superalgebra ${\mathbb K}e$ with null even part by means of $(D,x_0)$.\ \ \ \ \ $\square$ \\

If we consider $({\frak g},\omega)$ a nilpotent odd-symplectic Lie superalgebra, then the odd-symplectic Lie superalgebra $({\frak h},{\omega}_{\frak h})$ obtained in Proposition \ref{prop-inverse-impaire}  is also a nilpotent odd-symplectic Lie superalgebra. So we can give an inductive description of nilpotent odd-symplectic Lie superalgebras.\\

\begin{cor}
Every non-null nilpotent odd-symplectic Lie superalgebra is obtained from $\left\{0\right\}$ by a sequence of double extensions by the one-dimensional Lie algebra and/or generalized double extensions by the one-dimensional Lie superalgebra with null even part.\\
\end{cor}

\noindent\textbf{2.2 Inductive description of nilpotent even-symplectic Lie superalgebras}
\\

In this subsection we give inductive description of nilpotent even-symplectic Lie superalgebra after introducing the double extension and the generalized double extension of an even-symplectic Lie superalgebra by a one-dimensional Lie superalgebra.\\

\begin{prop}
Let $({\frak g},\omega)$ be an even-symplectic Lie superalgebra and $D$ a homogeneous superderivation on ${\frak g}$ of degree $\mid e\mid$. Define the bilinear form $\gamma: {\frak g}\times {\frak g}\longrightarrow {\mathbb K}$ for homogeneous elements by
\begin{equation*}
\gamma(X, Y) = - {(-1)}^{\mid e\mid}\Bigl[\omega\bigl(D(X), Y\bigl) + (-1)^{\mid e\mid \mid x\mid} \omega\bigl(X, D(Y)\bigl)\Bigl].
\end{equation*}
Then the ${\mathbb Z}_2$-graded vector space ${\mathbb K}e^* \oplus {\frak g}$ endowed with the following bracket
\begin{equation*}
 \quad [a e^*+X, b e^*+Y] :=  \gamma(X,Y) e^* + [X, Y]_{\frak g}, {\forall}_{(a e^*+ X),(b e^* +Y)\in {\mathbb K}e^*\oplus{\frak g}}
\end{equation*}is a Lie superalgebra. ${\mathbb K}e^*\oplus {\frak g}$ is called the central extension of ${\frak g}$ by ${\mathbb K}e^*$ by means of $\gamma$.\\
\end{prop}

\begin{thm}\label{g-d-e-even-symp}
Let $({\frak g},\omega)$ be an even-symplectic Lie superalgebra, $D\in {(Der({\frak g}))}_{\mid e\mid}$ and $x_0\in {\frak g}_{\bar{0}}$ such that
\begin{eqnarray}
0&=&(1 - {(-1)}^{\mid e\mid})D(x_0),\\
0&=&(1 - {(-1)}^{\mid e\mid})\omega(x_0,x_0),\\
0&=&(1 - {(-1)}^{\mid e\mid})(D^2- \frac{1}{2}\left[x_0,.\right]).
\end{eqnarray}
 Consider the following scalar 2-cocycle of ${\frak g}$ defined by
$$\theta(X,Y)= - {(-1)}^{\mid e\mid}\omega(\Bigl[D^*({(-1)}^{\mid e\mid \mid x\mid}D+ D^*) + (D + {(-1)}^{\mid e\mid (\mid x\mid + 1)} D^*)D\Bigl](X),Y), \ \ {\forall}_{X,Y\in {\frak g}},$$
where $D^*$ is the adjoint of $D$ with respect to $\omega$. If $\theta$ is a 2-coboundary of ${\frak g}$ (i.e there exist $b_0\in {\frak g}_{\bar{0}}$ such that $\theta(X,Y) = \omega(b_0, [X,Y]),\ {\forall}_{X,Y\in {\frak g}}$), which is the case if $D$ is an odd superderivaion  (we choose $b_0= \frac{1}{2}x_0$), then the linear map $\widetilde{D}$ on ${\mathbb K}e^*\oplus{\frak g}$ defined by
$$\widetilde{D}(a e^*+X):= D(X) + \omega(b_0, X)e^*,\ {\forall}_{(a e^*+X)\in {\mathbb K}e^*\oplus{\frak g}}$$
is a homogeneous superderivation on ${\mathbb K}e^*\oplus{\frak g}$ of degree $\mid e\mid$. Consequently, the ${\mathbb Z}_2$-graded vector space ${\frak t}= {\mathbb K}e^* \oplus {\frak g}\oplus {\mathbb K}e$ endowed with the following bracket defined by
\begin{eqnarray*}\label{produitimpaire}
[e, e] &=& \frac{1}{2}(1 - {(-1)}^{\mid e\mid})x_0 \ ; \quad [e, X] = \widetilde{D}(X);\\
 \quad [X, Y] &=& {\left[X,Y\right]}_{\frak g} +\gamma(X, Y) e^* ,\ {\forall}_{X, Y \in {\frak g}}.
\end{eqnarray*}
is a Lie superalgebra. Moreover the  skew-supersymmetric bilinear form $\Omega$ on ${\frak t}$ defined by
\begin{equation}
\Omega_{\vert_{\frak g\times \frak g}} = \omega,\ \ \Omega(e^*, e) = 1, \ \ \ \mbox{and}\ \ \Omega \ \mbox{ is null elsewhere},
\end{equation}
is an even-symplectic form on ${\frak t}$. $({\frak t}, \Omega)$ is an even-symplectic Lie superalgebra which we called the double extension (rep. generalized double extension )of $({\frak g}, \omega)$ by the one-dimensional Lie algebra (resp. the one-dimensional Lie superalgebra with null even part) by means of $(D,b_0)$ (resp. $(D,x_0)$).\\
\end{thm}

\begin{prop}\label{prop-inverse-paire}
Let $({\frak g}, \omega)$ be an even-symplectic Lie superalgebra different from the one-dimensional Lie superalgebra with null even part. If ${\frak z}({\frak g})\neq \left\{0\right\}$,
 then $({\frak g}, \omega)$ is a one of the three following propositions:
 \begin{enumerate}
 \item[(i)] Double extension of an even-symplectique Lie superalgebra by the one dimensional Lie algebra,
 \item[(ii)]Generalized double extension of  even-symplectic Lie superalgebra  by the one-dimensional Lie superalgebra with null even part,
 \item[(iii)] Orthogonal direct sum of a non-degenerate one-dimensional graded ideal of ${\frak g}$ with null even part and a non-degenerate graded ideal of ${\frak g}$ with null center.\\
\end{enumerate}
\end{prop}

\noindent\underline{\sl Proof} : Since ${\frak z}({\frak g})\neq \left\{0\right\}$, then we have to analyze two cases.\\

\textbf{First case:} We suppose that ${\frak z}({\frak g})\cap {\frak g}_{\bar{0}}\neq \left\{0\right\}$. We denote $I:= {\mathbb K}e^*$, where $e^*\in {\frak z}({\frak g})\cap {\frak g}_{\bar{0}}\setminus \left\{0\right\}$, and $J$ its orthogonal with respect to $\omega$. Since $\omega$ is even and skew-supersymmetric, it is clear that $I\subseteq J$. As $\omega$ is even and non-degenerate then there exists $e\in {\frak g}_{\bar{0}}$ such that $\omega(e^*,e)= 1$. From the non-degeneracy of $\omega$, there exists a graded vector space ${\frak h}$ of ${\frak g}$ such that $\omega(e^*, {\frak h}) = \{0\} = \omega(e, {\frak h}), J = I \oplus {\frak h}$ and ${\omega}_{\frak h} = \omega_{|_{\frak h\times\frak h}}$ is a non-degenerate form on ${\frak h}$. Since ${\frak h}$ is a graded vector space of ${\frak g}$ contained in the graded ideal $J$ then the bracket of ${\frak g}$ takes the form :
\begin{eqnarray*}
[X, Y] &=& \gamma(X, Y) e^* + \alpha(X,Y)\ \ \ \mbox{and}\ \ \  [e, X] = D(X) + \beta(X)e^*, {\forall}_{X,Y\in {\frak h}},
\end{eqnarray*}
where $\alpha(X,Y),D(X)\in {\frak h}$ and $\gamma(X,Y), \beta(X)\in {\mathbb K}$. Since  ${\frak g}$ is a Lie superalgebra, it comes directly that $({\frak h},\alpha)$ is a Lie superalgebra. Besides it not difficult to see that ${\omega}_{\frak h}$ is an even-symplectic from on $({\frak h},\alpha)$. Moreover, by using the fact that $\omega$ is a scalar 2-cocycle, we deduce that
$$\gamma(X, Y) = - \Bigl[{\omega}_{\frak h}\bigl(D(X), Y\bigl) + {\omega}_{\frak h}\bigl(X, D(Y)\bigl)\Bigl].$$
On the other hand, by the graded identity of Jacobi, it comes that $D$ is an even  superderivation of ${\frak h}$ and for homogeneous elements $X$ and $Y$ in ${\frak h}$, we have
\begin{equation}
\label{condition-cobord-deux}\beta(\left[X,Y\right]) = - {\omega}_{\frak h} (\Bigl[D^*(D + D^*) + (D + D^*)D\Bigl](X),Y),
\end{equation}
which prove that the bilinear form $\theta: {\frak h}\times {\frak h}\longrightarrow {\mathbb K}$ defined by
$$\theta(X,Y):= -{\omega}_{\frak h}(\Bigl[ D^*(D+ D^*) + (D + D^*)D\Bigl](X),Y), {\forall}_{X,Y \in {\frak h}},$$
is a scalar 2-coboundary of ${\frak h}$. Since ${\omega}_{\frak h}$ is non-degenerate and $\beta\in {\frak h}^*$, then there exists a unique $z_0\in {\frak h}$ such that $\beta= \omega(z_0,.)$. Now, applying Theorem \ref{g-d-e-even-symp}, we conclude that ${\frak g}$ is the double extension of $({\frak h}, {\omega}_{\frak h})$ by the one-dimensional Lie algebra ${\mathbb K}e$ by means of $(D,z_0)$.\\

\textbf{Second case:} We suppose that ${\frak z}({\frak g})\cap {\frak g}_{\bar{0}}= \left\{0\right\}$. Since ${\frak z}({\frak g})\neq \left\{0\right\}$, then we have ${\frak z}({\frak g})\cap {\frak g}_{\bar{1}}\neq  \left\{0\right\}$.

$\bullet$ If $dim ({\frak z}(\frak g)\cap {\frak g}_{\bar{1}})=1$, then ${\frak z}({\frak g})= {\frak z}({\frak g})\cap {\frak g}_{\bar{1}}={\mathbb K}e^*$. Set $\omega(e^*,e^*)= k$, where $k\in {\mathbb K}$. If $k\neq 0$, then we have $I:= {\mathbb K}e^*$ is a non-degenerate graded ideal of ${\frak g}$ and  ${\frak g}= I \oplus I^{\bot}$. By Lemma \ref{orthogonal-of-I-is-ideal}, it comes that $I^{\bot}$ is a graded ideal of ${\frak g}$. Since $I$ is non-degenerate and $dim ({\frak z}({\frak g}))=1$, we deduce that ${I}^{\bot}$ is  a non-degenerate graded ideal of ${\frak g}$ with null center. If now $k=0$, then we have $I:={\mathbb K}e^*$ is totally isotropic. As $\omega$ is even and non-degenerate then there exists $e\in {\frak g}_{\bar{1}}$ such that $\omega(e^*,e)= 1$. From the non-degeneracy of $\omega$, there exists a graded vector space ${\frak h}$ of ${\frak g}$ such that $\omega(e^*, {\frak h}) = \{0\} = \omega(e, {\frak h}), J = I \oplus {\frak h}$ and ${\omega}_{\frak h} = \omega_{|_{\frak h\times\frak h}}$ is a non-degenerate form on ${\frak h}$. Since ${\frak h}$ is a graded vector space of ${\frak g}$ contained in the graded ideal $J$ then the bracket of ${\frak g}$ takes the form :
\begin{eqnarray*}
[X, Y] &=& \gamma(X, Y) e^* + \alpha(X,Y), \ \ \ [e, X] = D(X) + \beta(X)e^*\ \ \mbox{and}\ \ \ \left[e,e\right]=x_0,
\end{eqnarray*}
for all $X,Y\in {\frak h}$, where $\alpha(X,Y),D(X)\in {\frak h}$, $x_0\in {\frak h}_{\bar{0}}$ and $\gamma(X,Y), \beta(X)\in {\mathbb K}$. Since  ${\frak g}$ is a Lie superalgebra, it comes directly that $({\frak h},\alpha)$ is a Lie superalgebra. Besides it not difficult to see that ${\omega}_{\frak h}$ is an even-symplectic from on ${\frak h}$. Moreover, by using the fact that $\omega$ is a scalar 2-cocycle, we deduce that $\beta(X)=\frac{1}{2}\omega(x_0,X)$ and
\begin{equation*}
\gamma(X, Y) = \Bigl[{\omega}_{\frak h}\bigl(D(X), Y\bigl) + (-1)^{\mid x\mid} {\omega}_{\frak h}\bigl(X, D(Y)\bigl)\Bigl], \ \forall X,Y\in {\frak h}.
\end{equation*}
On the other hand, by the graded identity of Jacobi, it comes that $\omega(x_0,x_0)=0$, $D(x_0)=0$, $D^2=\frac{1}{2}ad_{\frak h}(x_0)$ and $D$ is an odd superderivation of ${\frak h}$. Hence we conclude, by applying Theorem \ref{g-d-e-even-symp}, that $({\frak g}, \omega)$ is a generalized double extension of $({\frak h}, {\omega}_{\frak h})$ by the one-dimensional Lie superalgebra with null even part ${\mathbb K}e$ by means of $(D,x_0)$.

$\bullet$ If now we suppose that $dim ({\frak z}(\frak g)\cap {\frak g}_{\bar{1}})\geq 2$, then we can see easily that there exists $e^*\in {\frak z}({\frak g})\cap {\frak g}_{\bar{1}}\setminus \left\{0\right\}$ such that $\omega(e^*,e^*)=0$ i.e $I:= {\mathbb K}e^*$ is a totally isotropic graded ideal. Using the same reasoning when $dim({\frak z}(\frak g)\cap {\frak g}_{\bar{1}})=1$, we deduce that $({\frak g}, \omega)$ is a generalized double extension of an even-symplectic Lie superalgebra by the one-dimensional Lie superalgebra with null even part. \ \ \ \  $\square$ \\

\begin{cor}\label{nilpotent-case}
Every  nilpotent even-symplectic Lie superalgebra different from the one-dimensional Lie superalgebra with null even part is obtained as either $(i)$ or $(ii)$ of Proposition \ref{prop-inverse-paire}.\\
\end{cor}

\begin{rem}\label{h-is-nilpotent}
If $({\frak g},\omega)$ is a nilpotent even-symplectic Lie superalgebra, then the even-symplectic Lie superalegrba $({\frak h},{\omega}_{\frak h})$ obtained in Proposition \ref{prop-inverse-paire}  is also a nilpotent even-symplectic Lie superalegrba.\\
\end{rem}

Consequently, by the Corollary \ref{nilpotent-case} and the Remark \ref{h-is-nilpotent} we can give an inductive description of nilpotent even-symplectic Lie superalgebras.
Consider ${\frak U}$ the set composed by $\left\{0\right\}$ and the one-dimensional Lie superalgebra with null even part.\\

\begin{cor}\label{I-D-even-symp-nilpotent-Lie-SA}
Every nilpotent even-symplectic Lie superalgebra is either an element of ${\frak U}$ or obtained from a finite number of elements of  ${\frak U}$ by a finite sequence of double extensions by the one-dimensional Lie algebra and/or generalized double extensions by the one-dimensional Lie superalgebra with null even part.\\
\end{cor}

\noindent\underline{\sl Proof} : Let $({\frak g},\omega)$ be a nilpotent even-symplectic Lie superalgebra. To prove that ${\frak g}$ is obtained as the Corollary \ref{I-D-even-symp-nilpotent-Lie-SA}, we proceed by induction on the dimension $n$ of ${\frak g}$. If $n=1$, then we deduce, by the skew-symmetry of $\omega$, that ${\frak g}$ is the one-dimensional Lie superalgebra with null even part. Suppose that the Corollary is true for $dim {\frak g}< n$ with $n\geq 2$ and we consider $dim {\frak g}=n$. By the Corollary \ref{nilpotent-case} and the Remark \ref{h-is-nilpotent}, we obtain that ${\frak g}$ is a double extension or a generalized double extension of a nilpotent even-symplectic Lie superalgebra ${\frak h}$ by a one-dimensional Lie superalgebra. In this case $dim {\frak h}= dim {\frak g} - 2< n$. If ${\frak h}$ is neither $\left\{0\right\}$ nor the one-dimensional Lie superalgebra with null even part, then we apply the induction hypothesis to ${\frak h}$, which complete the proof.\ \ \ \  $\square$ \\


In the sequel,  we apply the same reasoning followed in~\cite{quadsymp} but not the same tool to give an inductive description of the homogeneous symplectic quadratic Lie superalgebras. In~\cite{quadsymp}, the authors used the quadratic generalized double extension of a quadratic  Lie superalgebra by a one-dimensional Lie superalgebra to give an inductive description of quadratic symplectic Lie superalgebras. However, in this work, we use the odd-quadratic generalized double extension of an odd-quadratic  Lie superalgebra by the one-dimensional Lie superalgebra with null even part to give an inductive description of the odd-quadratic odd-symplectic Lie superalgebras. Moreover, we introduce even (resp. odd)-quadratic generalized double extension of even (resp. odd)-quadratic Lie superalgebra by the abelian two-dimensional Lie superalgebra to give an inductive description of the even (resp. odd)-quadratic odd (resp. even)-symplectic Lie superalgebras.\\

The following lemma it will be very useful in the study of homogeneous quadratic symplectic Lie superalgebras.\\

\begin{lem}\label{center-non-null}
Let $({\frak g}, B,\omega)$ be a homogeneous quadratic symplectic Lie superalgebra, then ${\frak z}({\frak g}):= \left\{X\in {\frak g}, \left[X,{\frak g}\right]=\left\{0\right\}\right\}$ is non-null.\\
\end{lem}

\noindent\underline{\sl Proof}: Following Proposition \ref{caracterisation}, there exist $\Delta$ a homogeneous invertible skew-supersymmetric superderivation of ${\frak g}$ such that $\omega(X,Y):= B(\Delta(X),Y)$, $\forall X,Y\in {\frak g}$. We have to analyse two cases:\\

\textbf{First case}: Suppose that $\Delta$ is an even invertible superderivation on ${\frak g}$. Since the restriction ${\Delta\vert}_{{\frak g}_{\bar{0}}}$ of $\Delta$ to ${\frak g}_{\bar{0}}$ is an invertible derivation of ${\frak g}_{\bar{0}}$, then by \cite{jacobson1955}, it follows that ${\frak g}_{\bar{0}}$ is a nilpotent Lie algebra. Hence, by \cite{Kac1977}, ${\frak g}$ is a solvable Lie superalgebra. Meaning that $\left[{\frak g},{\frak g}\right]\neq {\frak g}$ and so ${\frak z}({\frak g})\neq \left\{0\right\}$ since ${\frak z}({\frak g})= {\left[{\frak g},{\frak g}\right]}^{\bot}$ with respect to $B$.\\

\textbf{Second case}: Suppose that $\Delta$ is an odd invertible superderivation on ${\frak g}$ and consider $\widetilde{\Delta}:= \left[\Delta,\Delta\right]$. Clearly that $\widetilde{\Delta}$ is an even invertible superderivation on ${\frak g}$ and consequently we have the result from the first case of this proof. \ \ \ \ $\square$ \\

\section{Odd-quadratic Lie superaglebras with homogeneous symplectic structures}

\noindent\textbf{3.1 Odd-quadratic odd-symplectic Lie superalgebras}
\\

The notion of generalized odd double extension of an odd-quadratic Lie superalgebra by the one-dimensional Lie superalgebra with null even part was introduced in \cite{albuquerque2009}. This notion will be very useful to give inductive description of odd-quadratic odd-symplectic Lie superalgebras. For, this reason, we start by recalling the generalized odd double extension of odd-quadratic Lie superalgebras introduced in \cite{albuquerque2009}. Let $({\frak g},B)$ be an odd-quadratic Lie superalgebra, ${\mathbb K} e$ the one-dimensional Lie superalgebra with null even part (i.e ${\mathbb K} e = {({\mathbb K} e)}_{\bar{1}}$), $k\in {\mathbb K}$, $x_0\in {\frak g}_{\bar{0}}$ and $\bar{D}$ an odd skew-symmetric superderivation on ${\frak g}$ such that:
\begin{equation}
\label{condition-quad-odd}\bar{D}(x_0)=0\ \ \mbox{and}\ \ {\bar{D}}^2= \frac{1}{2}ad_{\frak g}(x_0).
\end{equation}
On the ${\mathbb Z}_2$-graded vector space $\widetilde{\frak g}:= P({\mathbb K} e^*)\oplus {\frak g}\oplus {\mathbb K} e$, we define the following bracket
\begin{eqnarray}
\nonumber{\left[e,e\right]}_{\widetilde{\frak g}}&=& x_0 + k e^*,\\
\label{bracket-quad-odd}{\left[e,X\right]}_{\widetilde{\frak g}}&=& \bar{D}(X) + B(x_0,X)e^*,\\
\nonumber{\left[X,Y\right]}_{\widetilde{\frak g}}&=& {\left[X,Y\right]}_{\frak g} + B(\bar{D}(X),Y)e^*,
\end{eqnarray}
and the following symmetric bilinear form
\begin{equation}
\label{form-quad-odd}{\widetilde{B}\vert}_{{\frak g}\times {\frak g}}:= B, \ \ \widetilde{B}(e^*,e)=1 \ \ \mbox{and}\ \widetilde{B}\ \ \mbox{is null elsewhere.}
\end{equation}
Provided with bracket ${\left[,\right]}_{\widetilde{\frak g}}$ and the bilinear form $\widetilde{B}$, $\widetilde{\frak g}$ is an odd-quadratic Lie superalgebra. Following \cite{albuquerque2009}, $(\widetilde{{\frak g}}, \widetilde{B})$ is called the generalized odd double extension of $({\frak g},B)$ by ${\mathbb K} e$ by means of $(\bar{D},x_0,k)$.\\


Now, we consider the generalized odd double extension $(\widetilde{\frak g}, \widetilde{B})$ of the odd-quadratic Lie superalgebra $({\frak g},B)$ by the one-dimensional Lie superalgebra by means of $(\bar{D},x_0,k)$. We suppose that there exist an even invertible skew-supersymmetric superderivation $\delta$ on $({\frak g},B)$. By Proposition \ref{caracterisation}, $({\frak g},B,\omega)$ is an odd-quadratic odd-symplectic Lie superalgebra. Let $c_1\in {\frak g}_{\bar{1}}$ and $\lambda\in {\mathbb K}\setminus \left\{0\right\}$ such that

\begin{equation}\label{condition-oddquad-oddsymp}
B(c_1,x_0)=\lambda k ,\ \ 2\bar{D}(c_1)=\delta(x_0)+2\lambda x_0,\ \ \left[\delta, \bar{D}\right]+\lambda \bar{D}= ad_{\frak g}(c_1).
\end{equation}

Then, the following linear map $\Delta$ defined on $\widetilde{\frak g}$ by

\begin{equation}
\label{derivation-quad-odd}\Delta(e^*)= \lambda e^*,\ \ \Delta(X)=\delta(X) - B(c_1,X)e^*, \ \ \Delta(e)= -\lambda e + c_1,
\end{equation}

is an even invertible skew-symmetric superderivation on $(\widetilde{\frak g},\widetilde{B})$. Consequently, $(\widetilde{\frak g}, \widetilde{\omega})$, where $\widetilde{\omega}(.,.):= \widetilde{B}(\Delta(.),.)$, is an odd-symplectic Lie superalgebra. The odd-quadratic odd-symplectic Lie superalgebra $(\widetilde{\frak g}, \widetilde{B}, \widetilde{\omega})$ is called the generalized odd double extension of the odd-quadratique odd-symplectic Lie superalgebra $({\frak g}, B, \omega)$ by the one-dimensional Lie superalgebra by means of $(\bar{D},x_0,c_1,\lambda)$.\\

\begin{rem}
For an odd-quadratic odd-symplectic Lie superalgebra $({\frak g}, B, \omega)$, we have ${\frak z}({\frak g}_{\bar{0}})\neq \left\{0\right\}$. Moreover, following the parity and the invariance of $B$, it comes ${\frak z}({\frak g}_{\bar{0}})\subseteq {\frak z}({\frak g})\cap {\frak g}_{\bar{0}}$. So we have ${\frak z}({\frak g})\cap {\frak g}_{\bar{0}}\neq \left\{0\right\}$. For this reason, in this section, we consider only the generalized odd double extension of odd-quadratic odd-symplectic Lie superalgebra by the one-dimensional Lie superalgebra with null even part.\\
\end{rem}

\begin{prop}\label{inverse-oddquad-oddsymp}
Every odd-quadratic odd-symplectic Lie superalgebra  is a generalized odd double extension of an odd-quadratic odd-symplectic Lie superalgebra $({\frak h}, B_{\frak h}, {\omega}_{\frak h})$ by the one-dimensional Lie superalgebra with null even part.\\
\end{prop}
\noindent\underline{\sl Proof} : Let $({\frak g},B,\omega)$ be an odd-quadratic odd-symplectic Lie superalgebra and $\Delta$ the unique even invertible skew-supersymmetric superderivation of $({\frak g}, B)$ such that $\omega(X, Y) = B\bigl(\Delta(X), Y\bigl)$, for all $X, Y \in {\frak g}$. Since $\Delta$ is invertible, the field ${\mathbb K}$ is algebraically closed and $\Delta\bigl({\frak z}({\frak g})\cap {\frak g}_{\bar 0}\bigl) = {\frak z}({\frak g})\cap {\frak g}_{\bar 0}$, there exist a non zero scalar $\lambda$ and a non zero vector $e^*$ in ${\frak z}({\frak g})\cap {\frak g}_{\bar 0}$ such that $\Delta(e^*) = \lambda e^*$. If we denote by $I$ the graded ideal generated by $e^*$ and $J$ its orthogonal  with respect to $B$, then it is clear that $I\subseteq J$. Since $B$ is an odd and non-degenerate, then there exist $e\in {\frak g}_{\bar{1}}$ such that $B(e^*,e)=1$. The fact that $W:= {\mathbb K}e^*\oplus {\mathbb K}e$ is non-degenerate, yields that ${\frak g}= W\oplus {\frak h}$ where ${\frak h}= W^{\bot}$ with respect to $B$. It comes that $B_{\frak h}:= B_{{\frak h}\times {\frak h}}$ is non-degenerate. By the proof of Proposition 2.15 of  \cite{albuquerque2009}, $({\frak h}, {\left[,\right]}_{\frak h}:= p\circ {\left[,\right]\vert}_{{\frak h}\times {\frak h}}, B_{\frak h})$ is an odd-quadratic Lie superalgebra where $p: {\frak h}\oplus {\mathbb K}e^* \longrightarrow {\frak h}$ denotes the projection $p(X+ \alpha e^*)=X$, for all $X$ in ${\frak h}$. It comes also following the same proposition that ${\frak g}$ is the generalized odd double extension of $({\frak h}, {\left[,\right]}_{\frak h},B_{\frak h})$ by the one-dimensional Lie superalgebra with null even part ${\mathbb K}e$ by means of $(\bar{D},x_0,k)$, where $\bar{D}:= p\circ {ad_{\frak g}(e)\vert}_{\frak h}$ and $\left[e,e\right]= x_0 + k e^*$. Since ${\frak g}= P({\mathbb K}e^*)\oplus{\frak h}\oplus {\mathbb K}e$, $\Delta(J)\subseteq J$ and $J= P({\mathbb K}e^*)\oplus{\frak h}$ then we have:
$$\left\lbrace
\begin{array}{l}
\Delta(e^*)  = \lambda e^*,\\
\Delta(e)  = c_1 + \alpha e,\\
\Delta(X)  = \mu(X) e + \delta(X), \quad \forall X \in {\frak h},
\end{array}
\right.$$
 where $\alpha, \mu(X) \in {\mathbb K}$, $c_1\in {\frak h}_{\bar{1}}$ and $\delta(X) \in {\frak h}$. By the skew-supersymmetry of $\Delta$, we get $\delta$ is skew-symmetric with respect to $B_{\frak h}$, $\alpha = -\lambda$ and $\mu(X)= -B_{\frak h}(c_1,X), \ \forall X\in {\frak h}$.
In addition, the fact that $\Delta$ is an even invertible superderivation on ${\frak g}$, implies that
$$\Delta\left[X,Y\right]= \left[\Delta(X),Y\right] + \left[X,\Delta(Y)\right], \ \forall X,Y\in {\frak h}.$$
Consequently, we obtain $\delta$ is an even superderivation on ${\frak h}$ and $\left[\delta, D\right] + \lambda D= ad_{\frak h}(c_1)$. Moreover since $\Delta(\left[e,e\right])= 2\left[e, \Delta(e)\right]$, then we have $\delta(x_0) + 2\lambda x_0 = 2 D(c_1)$ and $B(x_0, c_1)=\lambda k$. Hence, we obtain that the conditions (\ref{condition-oddquad-oddsymp}) are satisfied. So, we deduce that $({\frak g}, B, \omega)$ is a generalized odd double extension of the odd-quadratic odd-symplectic Lie superalgebra $({\frak h}, B_{\frak h}, {\omega}_{\frak h})$, where ${\omega}_{\frak h} (X,Y)=B_{\frak h} \bigl(\delta(X), Y\bigl),\ \forall X, Y \in {\frak h}$, by the one-dimensional Lie superalgebra by means of $(\bar{D},x_0, c_1 \lambda)$.\ \ \ \ \ $\square$\\


It is clear that every 2-dimensional odd-quadratic odd-symplectic Lie superalgebra is abelian. Consequently, we have the following corollary.

\begin{cor}
Every non-null odd-quadratic odd-symplectic Lie superalgebra is obtained from $\left\{0\right\}$ by a finite sequence of generalized odd double extensions by the one-dimensional Lie superalgebra with null even part.\\
\end{cor}

\noindent\textbf{3.2 Odd-quadratic even-symplectic Lie superalgebras}
\\

In the following, we are going to give inductive description of odd-quadratic even-symplectic Lie superalgebras by using the notion of generalized double extension. To be done, we begin by introducing the notion of the generalized  double extension of odd-quadratic Lie superalgebras by the two-dimensional abelian Lie superalgebra.\\

\begin{thm}\label{theorem}
Let $({\frak g},B)$ be an odd-quadratic Lie superalgebra, $D$ and $\bar{D}$ two skew-supersymmetric superderivations on ${\frak g}$ which are respectively even and odd and let $x_0\in {\frak g}_{\bar{0}}$ and $x_1\in {\frak g}_{\bar{1}}$ such that the following conditions are satisfied
\begin{eqnarray}
\label{condition-odd}{\bar{D}}^2 &=& \frac{1}{2}ad_{\frak g}(x_0)\ \ \bar{D}(x_0)=0,\\
\label{condition-compatibility}\left[D,\bar{D}\right]&=& ad_{\frak g}(x_1),\ \ \  D(x_0)=2 \bar{D}(x_1), \ \ B(x_0,x_1) =0.
\end{eqnarray}
Then, the ${\mathbb Z}_{2}$-graded vector space ${\frak t}:= P({({\mathbb K}e_0\oplus {\mathbb K}e_1)}^*)\oplus {\frak g}\oplus {\mathbb K}e_0\oplus {\mathbb K}e_1$, where ${\mathbb K}e_0\oplus {\mathbb K}e_1$ is the two-dimensional abelian Lie superalgebra, endowed with the even skew-supersymmetric bilinear  map $\left[,\right]: {\frak t}\times {\frak t} \longrightarrow {\frak t}$ defined by
\begin{eqnarray*}
\label{extnesion-centrale}\left[X,Y\right]&:=& {\left[X,Y\right]}_{\frak g} + B(\bar{D}(X),Y)e_0^* + B(D(X),Y)e_1^*, \ \forall X,Y\in {\frak g},\\
\nonumber \left[e_0, X\right]&:=& D(X) - B(x_1,X){e_0}^*,\\
\nonumber \left[e_1,X\right]&:=& \bar{D}(X) + B(x_0,X){e_0}^* + B(x_1,X){e_1}^*,\\
\nonumber \left[e_0,e_1\right]&:=& x_1,\ \ \ \left[e_1,e_1\right]:= k {e_0}^* + x_0, \ k\in {\mathbb K},
\end{eqnarray*}
 is a Lie superalgebra. Moreover, the supersymmetric bilinear form $\gamma: {\frak t}\times {\frak t}\longrightarrow {\mathbb K}$ defined by:
\begin{equation}\label{structure-quadratique-impaire}
{\gamma\vert}_{{\frak g}\times {\frak g}}:= B, \ \ \gamma(e_0^*,e_0)=1= \gamma(e_1^*, e_1) \ \ \mbox{and}\ \gamma\ \ \mbox{is null elsewhere,}
\end{equation}
is an invariant odd scalar product ${\frak t}$. In this case, we say that $({\frak t},\gamma)$ is the generalized odd double extension of $({\frak g},B)$ by the two-dimensional abelian Lie superalgebra ${\mathbb K}e_0\oplus {\mathbb K}e_1$ by means of $(D, \bar{D}, x_0, x_1,k).$\\
\end{thm}

Following the hypotheses of the previous theorem, we remark that we can consider the generalized odd double extension of ${\frak g}$ by the one-dimensional Lie superalgebra with null even part ${\mathbb K}e_1$ by means of $(\bar{D},x_0)$. So, we can rewrite the previous theorem as follows.\\

\begin{prop}
Let $(\widetilde{\frak g}:= {\mathbb K}e_0^*\oplus {\frak g}\oplus {\mathbb K}e_1, {\left[,\right]}_{\widetilde{\frak g}}, \widetilde{B})$ be the generalized odd double extension of the odd-quadratic Lie superalgebra $({\frak g},B)$ by ${\mathbb K}e_1$ by means of $(\bar{D},x_0,k)$ (we suppose that $P({\mathbb K}e_1^*):= {\mathbb K}e_0^*$). Consider $D$ an even skew-symmetric superderivation on $({\frak g},B)$ and $x_1\in {\frak g}_{\bar{1}}$ such that:
\begin{equation}
\label{condition-prop}\left[D,\bar{D}\right]= ad_{\frak g}(x_1),\ \ \  D(x_0)=2 \bar{D}(x_1), \ \ B(x_0,x_1) =0.
\end{equation}
Then, the linear map $\widetilde{D}$ on $\widetilde{\frak g}$ defined by
$$ \left\lbrace
\begin{array}{l}
\widetilde{D}(X):= D(X) - B(x_1,X)e_0^*, \ \forall X\in {\frak g},\\
\widetilde{D}(e_1):=x_1,\\
\widetilde{D}(e_0^*):=0,
\end{array}
\right.$$
is an even skew-supersymmeric superderivation on $(\widetilde{\frak g},\widetilde{B})$. Consequently, the ${\mathbb Z}_2$-graded vector space ${\frak t}:= {\mathbb K}e_1^*\oplus \widetilde{\frak g}\oplus {\mathbb K}e_0$ (we suppose that $P({\mathbb K}e_0^*):={\mathbb K}e_1^* $) endowed with the skew-supersymmetric bilinear map $\left[,\right]:{\frak t}\times {\frak t}\longrightarrow {\frak t}$ defined by
\begin{eqnarray*}
\left[X,Y\right]&:=& {\left[X,Y\right]}_{\widetilde{\frak g}} + B(\widetilde{D}(X),Y)e_1^*,\\
\left[e_1,X\right]&:=& {\left[e_1,X\right]}_{\widetilde{\frak g}} + B(x_1,X)e_1^*,\\
\left[e_1,e_1\right]&:=& {\left[e_1,e_1\right]}_{\widetilde{\frak g}},\\
\left[e_0,X\right]&:=&\widetilde{D}(X),\ \ \left[e_0,e_1\right]:= \widetilde{D}(e_1),
\end{eqnarray*}
and the bilinear form $\gamma: {\frak t}\times {\frak t}\longrightarrow {\mathbb K}$ defined by
\begin{equation}
{\gamma\vert}_{\widetilde{\frak g}\times \widetilde{\frak g}}:= \widetilde{B},\ \ \gamma(e_0^*,e_0):=1 \ \ \mbox{and} \ \gamma\ \mbox{is null elsewhere}
\end{equation}
is an odd-quadratic Lie superalgebra. $({\frak t},\left[,\right],\gamma)$ is called the odd double extension of the odd-quadratic Lie superalgebra $(\widetilde{\frak g},\widetilde{B})$ by the Lie algebra ${\mathbb K}e_0$ by means $\widetilde{D}$.\\
\end{prop}

\noindent\underline{\sl Proof} : By a simple computation where we use conditions (\ref{condition-prop}), we deduce that $\widetilde{D}$ is an even skew-supersymmetric superderivation of $(\widetilde{\frak g},\widetilde{B})$. So, folowing \cite{albuquerque2009}, we can consider the double extension of $(\widetilde{\frak g},\widetilde{B})$ by the Lie algebra ${\mathbb K}e_0$ by means of $\widetilde{D}$.\ \ \ \ \ $\square$\\

 Now we are going to introduce the generalized odd double extension of an odd-quadratic even-symplectic Lie superalgebra by the two-dimensional abelian Lie superalgebra.\\

 \begin{prop}\label{g-d-e-quad-impaire-symp-paire}
 Let $({\frak g}, B, \omega)$ be an odd-quadratic even-symplectic Lie superalgebra and $\delta$ the unique odd invertible skew-symmetric superderivation on ${\frak g}$ such that $\omega(X,Y)=B(\delta(X),Y)$, $\forall X,Y\in {\frak g}$. Consider $({\frak t}, \gamma)$ the generalized odd double extension of $({\frak g}, B)$ by the two-dimensional abelian Lie superalgebra ${\mathbb K}e_0\oplus {\mathbb K}e_1$ (by means of $(D, \bar{D}, x_0, x_1,k)$). If there exist $c_0\in {\frak g}_{\bar{0}}$, $c_1\in {\frak g}_{\bar{1}}$, $\lambda\in {\mathbb K}\setminus \left\{0\right\}$ and $\alpha \in {\mathbb K}$ such that the following conditions are satisfied:
\begin{eqnarray}
\label{c1symp-paire}ad_{\frak g}(c_0)&:=& \left[\delta, \bar{D}\right] - \lambda D,\\
\label{c2symp-paire} ad_{\frak g}(c_1)&:=& \left[\delta, D\right] + \lambda \bar{D},\\
\label{c3symp-paire}\delta(x_0)&:=& -2 \bar{D}(c_0) +2 \lambda x_1,\\
\label{c4symp-paire}\delta(x_1)&:=& D(c_0) + \bar{D}(c_1) - \lambda x_0,\\
\label{c5symp-paire}\lambda k&:=& B(c_1, x_0) - 2 B(x_1, c_0),
\end{eqnarray}

Then the following linear map $\Delta$ on ${\frak t}$ defined by:

$$\Delta(e_0^*):= \lambda e_1^*,\ \  \Delta(e_1^*):= \lambda e_0^*,\ \Delta (e_0):= c_1 - \lambda e_1,\ \ \Delta(e_1):= \alpha e_0^* + c_0 + \lambda e_0,$$
$$ \Delta(X):= B(c_0,X)e_0^* + \delta(X) - B(c_1,X)e_1^*.$$

is an invertible odd skew-supersymmetric superderivation on $({\frak t},\gamma)$. Consequently, $({\frak t}, \gamma, \Omega)$ is an odd-quadratic even-symplectic Lie superalgebra, where $\Omega: {\frak t}\times {\frak t} \longrightarrow {\mathbb K}$ is defined by
$$\Omega(X,Y):=\gamma(\Delta(X),Y),\ \forall X,Y\in {\frak t}.$$\\
The odd-quadratic even-symplectic Lie superalgebra $({\frak t},\gamma,\Omega)$ is called the generalized odd double extension of $({\frak g},B,\omega)$ by the two-dimensional abelian Lie superalgebra ${\mathbb K}e_0\oplus {\mathbb K}e_1$ by means of $(D,\bar{D}, x_0, x_1, c_0,c_1,\lambda)$.\\
\end{prop}

\noindent\underline{\sl Proof}: Since $\delta$ is odd and invertible so is the map $\Delta$. Next, we have to see that
\begin{equation}
\label{Delta-derivation}\Delta\left[X,Y\right] = \left[\Delta(X),Y\right] + {(-1)}^{\mid x\mid} \left[X,\Delta(Y)\right], \ \forall X\in {\frak t}_{\mid x\mid},Y\in {\frak t}.
\end{equation}
Analyzing each term, we obtain  easily from conditions (\ref{c1symp-paire}) and (\ref{c2symp-paire})  and the fact that $B$ is invariant and $\delta$ is an odd skew-symmetric superderivation, that (\ref{Delta-derivation}) hold for $X\in {\frak g}_{\mid x\mid},Y\in {\frak g}$. In addition, from (\ref{c2symp-paire}) and (\ref{c4symp-paire}) (resp. (\ref{c1symp-paire}) and (\ref{c4symp-paire})), we deduce that (\ref{Delta-derivation}) hold for $e_0$ and $X\in{\frak g}$ (resp. $e_1$ and $X\in{\frak g}$)
Finally, using conditions (\ref{c4symp-paire}) and (\ref{c5symp-paire}) (resp. (\ref{c3symp-paire}) and (\ref{c5symp-paire})), we infer that (\ref{Delta-derivation}) hold for $e_0$ and $e_1$ (resp. $e_1$ and $e_1$).
To complete the proof, we have to ensure that the odd superderivation $\Delta$ is skew-symmetric. As $\delta$ is an odd skew-symmetric superderivation of $({\frak g},B)$, it follows that $B(\Delta(X),Y) = -{(-1)}^{\mid x\mid}B(X,\Delta(Y)), \ \forall X\in {\frak g}_{\mid x\mid},Y\in{\frak g}$. The other cases come trivially, so the proof is complete.\ \ \ \ \ $\square$ \\

Our next goal is to establish the converse of Proposition \ref{g-d-e-quad-impaire-symp-paire}.

\begin{prop}\label{inverse-odd-quad-even-symp-Lie}
Every odd-quadratic even-symplectic Lie superalgebra $({\frak g}, B, \omega)$ is a generalized odd double extension of an odd-quadratic even-symplectic Lie superalgerba $({\frak h}, B_{\frak h}, {\omega}_{\frak h})$, where $dim {\frak h}= dim {\frak g}-4$, by the two-dimensional abelian Lie superalgebra.\\
\end{prop}

\noindent\underline{\sl Proof} : Let $\Delta$ be the unique odd invertible skew-symmetric superderivation on ${\frak g}$ such that $\omega(X,Y)= B(\Delta(X),Y),\ \forall X,Y\in {\frak g}$. The fact that the field ${\mathbb K}$ is algebraically closed, $\Delta$ is inversible and $\Delta ({\frak z}({\frak g}))= {\frak z}({\frak g})$, implies that there exists $e^*\in {\frak z}({\frak g})\setminus \left\{0\right\}$ such that $\Delta (e^*)= \lambda e^*$, where $\lambda\in {\mathbb K}\setminus\left\{0\right\}$. More precisely, we have $e^*=e_0^*+e_1^*$, where $e_0^*\in {\frak z}({\frak g})\cap {\frak g}_{\bar{0}}\setminus \left\{0\right\}$ and $e_1^*\in {\frak z}({\frak g})\cap {\frak g}_{\bar{1}}\setminus \left\{0\right\}$, such that $\Delta(e_0^*)= \lambda e_1^*$ and $\Delta(e_1^*)= \lambda e_0^*$. We denote the graded ideal $I:={\mathbb K}e_0^*\oplus {\mathbb K}e_1^*$ and $J$ its orthogonal with respect to $B$. Since $\omega$ is skew-supersymmetric, we obtain  $ B(e_0^*,e_1^*)= 0$ and so we have $I\subseteq J$. As $B$ is odd and  non-degenerate, then there exist $e_0\in {\frak g}_{\bar{0}}$ and $e_1\in {\frak g}_{\bar{1}}$ such that  $B(e_0^*, e_1)\neq 0$ and $B(e_1^*,e_0)\neq 0$. We can choose $e_0$ and $e_1$ such that
 $$B(e_0, e_1^*)= 1 = B(e_1, e_0^*) \ \mbox{and}\ B(e_0,e_1)=0.$$
 Consider the four-dimensional ${\mathbb Z}_2$-graded vector subspace $W:= I\oplus {\mathbb K}e_0\oplus {\mathbb K}e_1$ of ${\frak g}$. Since ${B\vert}_{W\times W}$ is not degenerate, we have ${\frak g}= W\oplus {\frak h}$, where ${\frak h}$ is the orthogonal of $W$ with respect to $B$. It comes that $B_{\frak h}= {B\vert}_{{\frak h}\times {\frak h}}$ is non-degenerate and $J=I\oplus {\frak h}$. As ${\frak h}$ is a graded vector subspace of ${\frak g}$ contained in $J$, we have  that
 $$\left[X,Y\right]= \alpha(X,Y) + \nu(X,Y)e_0^* + \nu'(X,Y)e_1^*, \ \forall\, X,Y\in {\frak h}$$
where $\alpha(X,Y)\in {\frak h}$ and $\nu(X,Y), \nu'(X,Y)\in {\mathbb K}$. Since $({\frak g},\left[,\right])$ is a Lie superalgebra, we deduce that $({\frak h}, \alpha)$ is a Lie superalgebra. Besides, by the invariance of $B$, it comes directly that
$$B_{\frak h}(\alpha(X,Y),Z)= B_{\frak h}(X,\alpha(Y,Z)), \ \forall X,Y,Z\in {\frak h}.$$
Then, we conclude that $({\frak h}, \alpha,B_{\frak h})$ is an odd-quadratic Lie superalgebra. For $X\in {\frak h}$, we have
 $$\left[e_0,X\right]= D(X) + f(X)e_0^* + h(X)e_1^*\ \ \mbox{and}\ \ \left[e_1,X\right]= \bar{D}(X) + F(X)e_0^* + H(X)e_1^*,$$
where $D(X), \bar{D}(X)\in {\frak h}$ and $f(X), h(X), F(X), H(X)\in {\mathbb K}$. Moreover, by the invariance of $B$, we obtain  $B(e_0, \left[e_0,e_1\right])=0= B(e_0^*, \left[e_0,e_1\right])$ and it follows that $\left[e_0,e_1\right]= x_1$, where $x_1\in {\frak h}_{\bar{1}}$. Once more, by the invariance of $B$, we get $B(e_1^*, \left[e_1,e_1\right])=0$ and  consequently $\left[e_1,e_1\right]= ke_0^* + x_0$, where $x_0\in {\frak h}_{\bar{0}}$ and $k\in {\mathbb K}$.\\

\textbf{Claim} $D$ and $\bar{D}$ are two skew-supersymmetric superderivations of $({\frak h},B_{\frak h})$ such that the conditions (\ref{condition-odd}) and (\ref{condition-compatibility}) of Theorem \ref{theorem} are satisfied.\\

\noindent\underline{\sl Proof of the Claim} : It is clear that $D$ and $\bar{D}$ are two endomorphisms which are respectively even and odd. Moreover, for $e_i$ and $X,Y\in {\frak h}$, where $i\in \left\{0,1\right\}$, it comes by the graded Jacobi identity  and the invariance of $B$ that $D$ and $\bar{D}$  are two skew-symmetric superderivations of ${\frak h}$ such that the conditions (\ref{condition-odd}) are satisfied. In addition the conditions (\ref{condition-compatibility}) are obtained by the non-degeneracy of $B$ and the graded Jacobi identity for $e_0, e_1, x\in {\frak h}$ and $e_0,e_0,e_1$. Consequently, by Theorem \ref{theorem}, we obtain that $({\frak g}, B)$ is a generalized odd double extension of the odd-quadratic  Lie superalgebra $({\frak h}, B_{\frak h})$ by the two-dimensional abelian Lie superalgebra ${\mathbb K}e_0\oplus {\mathbb K}e_1$.\\

 Now, the fact that $\Delta$ is skew-symmetric and $\Delta(I) \subset I$, implies that $\Delta(J) \subset J$. Consequently, we have
\begin{eqnarray*}
\Delta(e_0)&=& i e_1^* + c_1 + je_1\\
\Delta(e_1)&=& le_0^* + c_0 + m e_0\\
\Delta(X) &=& \mu(X)e_0^* + \delta(X) + \mu'(X)e_1^*,
\end{eqnarray*}
where $\delta(X)\in {\frak h}$, $c_0\in {\frak h}_{\bar{0}}$, $c_1\in {\frak h}_{\bar{1}}$ $i,j,l,m,\mu(X),\mu'(X)\in {\mathbb K}$. The fact that $\Delta$ is skew-symmetric implies that $\delta$ is skew-symmetric with respect to $B_{\frak h}$, $i=0$, $j= -\lambda$, $m= \lambda$, $\mu(X)= B_{\frak h}(c_0,X)$ and $\mu'(X)= - B_{\frak h}(c_1,X)$. Moreover, since $\Delta$ is an odd invertible superderivation, it comes that $\delta$ is an odd invertible superderivation of ${\frak h}$ and the conditions (\ref{c1symp-paire})-(\ref{c5symp-paire}) are satisfied. So, we deduce that $({\frak h}, B_{\frak h})$ provided with the bilinear ${\omega}_{\frak h}:{\frak h}\times {\frak h}\longrightarrow {\mathbb K}$ defined by $\omega_{\frak h}(X,Y)= B_{\frak h}(\delta(X),Y)$, $\forall X,Y\in {\frak h}$ is an odd-quadratic even-symplectic Lie superalgebra and $({\frak g}, B, \omega)$ is the generalized odd double extension of $({\frak h},{B}_{\frak h},{\omega}_{\frak h})$ by means of $( D,\bar{D}, x_0,x_1,c_0,c_1, \lambda)$.\ \ \ \ $\square$ \\

 \begin{rem}
Since for any odd-quadratic even-symplectic Lie superalgebra $({\frak g},B,\omega)$ we must have $dim {\frak g}_{\bar{0}}= dim {\frak g}_{\bar{1}}$ and $dim {\frak g}_{\bar{0}}$ is even, so we deduce that even-quadratic odd-symplectic Lie superalgebras exist only starting from dimension four. Moreover, we prove that all 4-dimensional even-quadratic odd-symplectic Lie superalgebra are obtained from $\left\{0\right\}$ by generalized double extension by the 2-dimensional abelian Lie superalgebra. Indeed, let us consider $({\frak g}:= {\mathbb K}e_0\oplus {\mathbb K}f_0\oplus {\mathbb K}e_1\oplus {\mathbb K}f_1,B,\omega)$ an odd-quadratic even-symplectic Lie superalgebras and $\Delta$ the unique invertible skew-super-symmetric superderivation of ${\frak g}$ such that $\omega(x,y)=B(\Delta(x),y)$, $\forall x,y\in {\frak g}$. Following the Lemma \ref{center-non-null}, we can suppose that $e= e_0+ e_1\in {\frak z}({\frak g})$ such that $\Delta(e_0)=\lambda e_1$ and $\Delta(e_1)=\lambda e_0$ where $\lambda$ is a non-zero scalar. So we deduce easily that $I:={\mathbb K}e_0\oplus {\mathbb K}e_1$ is totally isotropic with respect to $\omega$ and with respect to $B$. On the other hand suppose that $\left[f_1,f_1\right]= \alpha e_0 + \beta f_0$ and $\left[f_0,f_1\right]= \overline{\alpha}e_1 + \overline{\beta}f_1$, where $\alpha,\beta,\overline{\alpha},\overline{\beta}\in {\mathbb K}$. Using the fact that $I$ is totally isotropic with respect to $\omega$ and $\omega(e_0, \left[f_1,f_1\right])=0$, $\omega(e_1, \left[f_0,f_1\right])=0$ we deduce $\beta=0$ and $\overline{\beta}=0$. In addition since $I$ is totally isotropic with respect to $B$ and $B(f_0,\left[f_0,f_1\right])=0$, it comes $\overline{\alpha}=0$. So we deduce that ${\frak g}$ is obtained for $\left\{0\right\}$ by generalized double extension by ${\mathbb K}f_0\oplus {\mathbb K}f_1$ by means of $(0,0,0,0,0,0,\alpha)$.\\
\end{rem}

\begin{cor}
Let $({\frak g},B,\omega)$ be a non-null odd-quadratic even-symplectic Lie superalgebra. Then, ${\frak g}$ is obtained from $\left\{0\right\}$ by a finite sequence of generalized odd double extension by the two-dimensional abelian Lie superalgebra.\\
\end{cor}

\section{Even-quadratic odd-symplectic Lie superalgebras}

To give inductive description of even-quadratic odd-symplectic Lie superalgebras, we will formalize the notion of generalized double extension of an even-quadratic odd-symplectic Lie superalgebra by the two-dimensional abelian Lie superalgebra. We start by introducing the generalized double extension of an even-quadratic Lie superalgebra by the two-dimensional abelian Lie superalgebra.\\

\begin{thm}\label{theorem-even}
Let $({\frak g},B)$ be an even-quadratic Lie superalgebra, $D$ and $\bar{D}$ two skew-symmetric superderivations on ${\frak g}$ which are respectively even and odd and let $x_0\in {\frak g}_{\bar{0}}$ and $x_1\in {\frak g}_{\bar{1}}$ such that the following conditions are satisfied
\begin{eqnarray}
\label{condition-even}{\bar{D}}^2 &=& \frac{1}{2}ad_{\frak g}(x_0), \ \ \ \bar{D}(x_0)=0,\ \ B(x_0,x_0) =0,\\
\label{condition-compatibility-even}\left[D,\bar{D}\right]&=& ad_{\frak g}(x_1),\ \ \  D(x_0)=2 \bar{D}(x_1).
\end{eqnarray}
Then, the ${\mathbb Z}_{2}$-graded vector space ${\frak t}:= {\mathbb K}e_0^*\oplus {\mathbb K}e_1^*\oplus {\frak g}\oplus {\mathbb K}e_0\oplus {\mathbb K}e_1$, where ${\mathbb K}e_0\oplus {\mathbb K}e_1$ is the two-dimensional abelian Lie superalgebra, endowed with the even skew-symmetric bilinear map $\left[,\right]: {\frak t}\times {\frak t} \longrightarrow {\frak t}$ defined by
\begin{eqnarray*}
\label{extnesion-centrale}\left[X,Y\right]&:=& {\left[X,Y\right]}_{\frak g} - B(\bar{D}(X),Y)e_1^* + B(D(X),Y)e_0^*, \ \forall X,Y\in {\frak g},\\
\nonumber \left[e_0, X\right]&:=& D(X) + B(x_1,X){e_1}^*,\\
\nonumber \left[e_1,X\right]&:=& \bar{D}(X) - B(x_0,X){e_1}^* + B(x_1,X){e_0}^*,\\
\nonumber \left[e_1,e_1\right]&:=& k {e_0}^* + x_0, \ \ \left[e_0,e_1\right]:= k e_1^* + x_1, \ k\in {\mathbb K},
\end{eqnarray*}
is a Lie superalgebra. Moreover, the supersymmetric bilinear form $\gamma: {\frak t}\times {\frak t}\longrightarrow {\mathbb K}$ defined by:
\begin{equation}\label{structure-quadratique-impaire}
{\gamma\vert}_{{\frak g}\times {\frak g}}:= B, \ \ \gamma(e_0^*,e_0)=1= \gamma(e_1^*, e_1) \ \ \mbox{and}\ \gamma\ \ \mbox{is null elsewhere,}
\end{equation}
is an invariant scalar product ${\frak t}$. In this case, we say that $({\frak t},\gamma)$ is the generalized double extension of $({\frak g},B)$ by the two-dimensional abelian Lie superalgebra ${\mathbb K}e_0\oplus {\mathbb K}e_1$ by means of $(D, \bar{D}, x_0, x_1,k).$\\
\end{thm}

\begin{example}
Consider the Lie superalgebra ${\frak L}$ introduced in Example \ref{example-nilpotent} and ${\frak g}= {\frak L}^*\oplus {\frak L}$ the trivial double extension of $\left\{0\right\}$ by  ${\frak L}$. On ${\frak L}$ we consider the even superderivation $D$ defined in Example \ref{example-nilpotent} and we define the odd linear maps $\bar{D}$ by
$$\bar{D}(l_0)=0,\ \ \bar{D}(k_0)=2k_1, \ \ \bar{D}(l_1)= l_0,\ \ \bar{D}k_1=0.$$
Doing some calculations, we can see easily that $\bar{D}$ is an odd superderivation on ${\frak L}$ such that $\bar{D}^2=0$ and $\left[D,\bar{D}\right]= \left\{0\right\}$. Now applying Lemma \ref{lemme-du-socle}, we obtain that $\widetilde{D}$ and $\widetilde{\bar{D}}$ are two superderivations on ${\frak g}$ which are respectively even and odd such that $\widetilde{\bar{D}}^2=0$ and $\left[\widetilde{D},\widetilde{\bar{D}}\right]=\left\{0\right\}$. Consequently, by Theorem \ref{theorem-even}, we can consider the generalized double extension of ${\frak g}$ by the two-dimensional abelian Lie superalgebra by means of $(\widetilde{D}, \widetilde{\bar{D}},0,0,0)$.\\
\end{example}

Now we are going to introduce the generalized double extension of an even-quadratic odd-symplectic Lie superalgebra by the two-dimensional abelian Lie superalgebra.\\

 \begin{prop}\label{g-d-e-quad-paire-symp-impaire}
 Let $({\frak g},B,\omega)$ be an even-quadratic odd-symplectic Lie superalgebra and $\delta$ the unique odd invertible skew-symmetric superderivation on ${\frak g}$ such that $\omega(X,Y)=B(\delta(X),Y)$, $\forall X,Y\in {\frak g}$. Consider $({\frak t}, \gamma)$ the generalized double extension of $({\frak g}, B)$ by ${\mathbb K}e_0\oplus {\mathbb K}e_1$ by means of $(D, \bar{D}, x_0, x_1,k)$. If there exist $c_0\in {\frak g}_{\bar{0}}$, $c_1\in {\frak g}_{\bar{1}}$, $\lambda\in {\mathbb K}\setminus \left\{0\right\}$ and $\alpha \in {\mathbb K}$ such that the following conditions are satisfied:
\begin{eqnarray}
\label{c1symp-impaire}\left[\delta, D\right] - \lambda \bar{D} &=& ad_{\frak a}(c_1),\\
\label{c2symp-impaire}\left[\delta, \bar{D}\right] + \lambda D &=& ad_{\frak a}(c_0), \\
\label{c3symp-impaire}\delta(x_0) &=& - 2 \lambda x_1 - 2 \bar{D}(c_0), \\
\label{c4symp-impaire}\delta(x_1) &=&  \lambda x_0 + D(c_0) + \bar{D}(c_1),\\
\label{c5symp-impaire} \lambda k &=& B(c_0,x_0).
\end{eqnarray}
Then the following linear map $\Delta$ defined on ${\frak t}$ by:

$$\Delta(e_0^*):= \lambda e_1^*,\ \  \Delta(e_1^*):= \lambda e_0^*,\ \Delta (e_0):=\alpha e_1^* +c_1 + \lambda e_1,\ \ \Delta(e_1):= -\alpha e_0^* + c_0 - \lambda e_0,$$
$$\Delta(X):= -B(c_1,X)e_0^* + \delta(X) -B(c_0,X)e_1^*,$$

is an invertible odd skew-symmetric superderivation on $({\frak t},\gamma)$. Consequently, $({\frak t}, \gamma, \Omega)$ is an even-quadratic odd-symplectic Lie superalgebra, where $\Omega: {\frak t}\times {\frak t} \longrightarrow {\mathbb K}$ is defined by
$$\Omega(X,Y):=B(\Delta(X),Y),\ \forall X,Y\in {\frak t}.$$\\
The even-quadratic odd-symplectic Lie superalgebra $({\frak t},\gamma,\Omega)$ is called the generalized double extension of $({\frak g}, B, \omega)$ by the two-dimensional abelian Lie superalgebra ${\mathbb K}e_0\oplus {\mathbb K}e_1$ by means of $(D,\bar{D}, x_0, x_1, c_0,c_1, \lambda)$.\\
\end{prop}

Our next goal is to establish the converse of Proposition \ref{g-d-e-quad-paire-symp-impaire}.

\begin{prop}\label{inverse-even-quad-odd-symp-Lie}
Every even-quadratic odd-symplectic Lie superalgebra $({\frak g}, B,\omega)$ is a generalized double extension of an even-quadratic odd-symplectic Lie superalgebra $({\frak h},B_{\frak h}, {\omega}_{\frak h})$ by the two-dimensional abelian Lie superalgebra.\\
\end{prop}

\noindent\underline{\sl Proof} : Let $\Delta$ be the unique odd invertible skew-symmetric superderivation on ${\frak g}$ such that $\omega(X,Y)= B(\Delta(X),Y),\ \forall X,Y\in {\frak g}$. The fact that the field ${\mathbb K}$ is algebraically closed, $\Delta$ is invertible and $\Delta ({\frak z}({\frak g}))= {\frak z}({\frak g})$, implies that there exists $e^*\in {\frak z}({\frak g})\setminus \left\{0\right\}$ such that $\Delta (e^*)= \lambda e^*$, where $\lambda\in {\mathbb K}\setminus\left\{0\right\}$. More precisely, we have $e^*=e_0^*+e_1^*$, where $e_0^*\in {\frak z}({\frak g})\cap {\frak g}_{\bar{0}}\setminus \left\{0\right\}$ and $e_1^*\in {\frak z}({\frak g})\cap {\frak g}_{\bar{1}}\setminus \left\{0\right\}$, such that $\Delta(e_0^*)= \lambda e_1^*$ and $\Delta(e_1^*)= \lambda e_0^*$. We denote the graded ideal $I:={\mathbb K}e_0^*\oplus {\mathbb K}e_1^*$. Since,
$$\lambda B(e_0^*,e_0^*)= \omega(e_1^*,e_0^*)= -\omega(e_0^*,e_1^*)= -\lambda B(e_1^*,e_1^*)=0,$$ it follows that $I$ is totally isotropic. Consequently  $I\subseteq J$, where $J$ is the orthogonal of $I$ with respect to $B$. Now,  as $B$ is even and non-degenerate, there exists $e_0\in {\frak g}_{\bar{0}}$ such that $B(e_0,e_0)=0$ and $B(e_0^*,e_0)=1$ and there exists $e_1\in {\frak g}_{\bar{1}}$ such that  $B(e_1^*, e_1)= 1$. Since the ${\mathbb Z}_2$-graded vector subspace $W:= I\oplus {\mathbb K}e_0\oplus {\mathbb K}e_1$ of ${\frak g}$ is not degenerate, it follows that ${\frak g}= W\oplus {\frak h}$, where ${\frak h}$ is the orthogonal of $W$ with respect to $B$. It comes that $B_{\frak h}= {B\vert}_{{\frak h}\times {\frak h}}$ is non-degenerate. Reasoning similar to Proposition \ref{inverse-odd-quad-even-symp-Lie}, we infer that $({\frak h}, {\left[,\right]\vert}_{{\frak h}\times {\frak h}},B_{\frak h},{\omega}_{\frak h})$ is an even-quadratic odd-symplectic Lie superalgebra such that ${\omega}_{\frak h}: {\frak h}\times {\frak h}\longrightarrow {\mathbb K}$ defined by:
$${\omega}_{\frak h}(X,Y):= B_{\frak h}(\delta(X),Y), \forall X,Y\in {\frak h},$$
where $\delta:= p_1\circ {\Delta\vert}_{\frak h}$ and $p_1: {\frak h}\oplus {\mathbb K}e_0^*\oplus {\mathbb K}e_1^* \longrightarrow {\frak h}$ denotes the projection $p_1(X+ \alpha e_0^*+ \beta e_1^*)= X,$ for all $X\in {\frak h}$ and $\alpha,\beta\in {\mathbb K}$. Once more reasoning similar to Proposition \ref{inverse-odd-quad-even-symp-Lie}, we obtain that ${\frak g}$ is the generalized double extension of $({\frak h}, {\left[,\right]}_{\frak h}:={\left[,\right]\vert}_{{\frak h}\times {\frak h}},B_{\frak h},{\omega}_{\frak h})$  by the two-dimensional abelian Lie superalgebra ${\mathbb K}e_0\oplus {\mathbb K}e_1$ by means $(D,\bar{D},x_0,x_1,c_0,c_1,\lambda)$, where $D= p_1\circ {ad_{\frak g}(e_0)\vert}_{\frak h}$, $\bar{D}=p_1\circ {ad_{\frak g}(e_1)\vert}_{\frak h}$, $\left[e_1,e_1\right]=x_0 + \frac{B_{\frak h}(x_0,c_0)}{\lambda}e_0^*$, $\left[e_0,e_1\right]=x_1 + \frac{B_{\frak h}(x_0,c_0)}{\lambda}e_1^*$, where $x_0,c_0\in {\frak h}_{\bar{0}}$ and $x_1,c_1\in {\frak h}_{\bar{1}}$, $c_0= p_2\circ \Delta(e_1)$ and $c_1= p_2\circ \Delta(e_0)$, where $p_2: {\frak g}\longrightarrow {\frak h}$ denotes the projection $p_2(\alpha e_0^* + \beta e_1^*+X+\alpha'e_0 + \beta'e_1)= X$, for all $X\in {\frak h}$ and $\alpha, \beta,\alpha',\beta'\in {\mathbb K}$. \ \ \ \ $\square$\\

\begin{rem}
Since for any even-quadratic odd-symplectic Lie superalgebra $({\frak g},B,\omega)$ we must have $dim {\frak g}_{\bar{0}}= dim {\frak g}_{\bar{1}}$ and $dim {\frak g}_{\bar{1}}$ is even, so we deduce that even-quadratic odd-symplectic Lie superalgebras exist only starting from dimension four. Moreover, we prove that all 4-dimensional even-quadratic odd-symplectic Lie superalgebra are obtained from $\left\{0\right\}$ by the generalized double extension by the 2-dimensional abelian Lie super-algebra. Indeed, let us consider that $({\frak g}:= {\mathbb K}e_0\oplus {\mathbb K}f_0\oplus {\mathbb K}e_1\oplus {\mathbb K}f_1,B,\omega)$ is an even-quadratic odd-symplectic Lie superalgebras and $\Delta$ the unique invertible skew-super-symmetric superderivation of ${\frak g}$ such that $\omega(x,y)=B(\Delta(x),y)$, $\forall x,y\in {\frak g}$. Following the Lemma \ref{center-non-null}, we can suppose that $e= e_0+ e_1\in {\frak z}({\frak g})$ such that $\Delta(e_0)=\lambda e_1$ and $\Delta(e_1)=\lambda e_0$ where $\lambda$ is a non-nul scalar. So, we deduce that $I:={\mathbb K}e_0\oplus {\mathbb K}e_1$ is totally isotropic with respect to $B$. On the other hand, we suppose that $\left[f_1,f_1\right]= \alpha e_0 + \beta f_0$ and $\left[f_0,f_1\right]= \overline{\alpha}e_1 + \overline{\beta}f_1$, where $\alpha,\beta,\overline{\alpha},\overline{\beta}\in {\mathbb K}$. Using the fact that $I$ is totally isotropic, $B(e_1, \left[f_0,f_1\right])=0$, $B(e_0, \left[f_1,f_1\right])=0$ and $B(f_0, \left[f_1,f_1\right])= B(\left[f_0,f_1\right],f_1)$, it comes that  $\overline{\beta}=0$, $\beta=0$ and $\alpha= \overline{\alpha}$. So we deduce that ${\frak g}$ is obtained for $\left\{0\right\}$ by that generalized double extension by ${\mathbb K}f_0\oplus {\mathbb K}f_1$ by means of $(0,0,0,0,,0,0,\alpha)$.\\
\end{rem}

\begin{cor}
Let $({\frak g},B,\omega)$ be a non-null even-quadratic odd-symplectic Lie superalgebra. Then, ${\frak g}$ is obtained from $\left\{0\right\}$ by a finite sequence of generalized double extension by the two-dimensional abelian Lie superalgebra.\\
\end{cor}

\section{Homogeneous-quadratic homogeneous-symplectic Lie superalgebras and homogeneous Manin superalgebras}

In \cite{quadsymp}, Theorem 5.7, it was proved that every even-quadratic even-symplectic Lie superalgebra over an algebraically closed field ${\mathbb K}$ is a special even-symplectic even-Manin superalgebra. In the same paper, the concept of generalized double extension of a special even-symplectic even-Manin superalgebra was introduced and it was used to describe even-quadratic even-symplectic Lie superalgebras. In this section, we will proceed similar to \cite{quadsymp}, to give inductive description of homogeneous-quadratic homogeneous-symplectic Lie superalgebras by using homogeneous-Manin superalgebras.\\

\begin{defn}
A homogeneous (i.e even or odd)-Manin superalgebra is a Lie superalgebra ${\frak g}$ which satisfy the following two conditions:
\begin{enumerate}
\item[(i)] ${\frak g}= {\frak a}\oplus {\frak b}$, where ${\frak a}$ and ${\frak b}$ are two Lie sub-superalgebras of ${\frak g}$ and ${\frak g}$ is the direct sum of the vector subspaces ${\frak a}$ and ${\frak b}$.
\item[(ii)] there exist a homogeneous (i.e even or odd)-quadratic structure $B$ on ${\frak g}$ such that ${\frak a}$ and ${\frak b}$ are totally isotropic with respect to $B$.\\
\end{enumerate}
\end{defn}

\begin{defn}
A special homogeneous-symplectic homogeneous-Manin superalgebra  is a homogeneous-Manin superalgebra $({\frak g}= {\frak a}\oplus {\frak b},B)$ provided with a homogeneous-symplectic structure $\omega$ such that $\omega({\frak a},{\frak a})= \omega({\frak b},{\frak b})= \left\{0\right\}$.\\
\end{defn}

The following lemma can be proved in the same way of Lemma 5.6 \cite{quadsymp}.\\

\begin{lem}
A homogeneous-symplectic homogeneous-Manin superalgebra $({\frak g}={\frak a}\oplus {\frak b},B, \omega)$ is a special homogeneous-symplectic homogeneous-Manin superalgebra if and only if the homogeneous invertible skew-supersymmetric superderivation $\Delta$ such that $\omega(X,Y)= B(\Delta(X),Y),\ \forall X,Y\in {\frak g}$ satisfies $\Delta({\frak a})\subseteq {\frak a}$ and $\Delta({\frak b})\subseteq {\frak b}$.\\
\end{lem}

\begin{thm}\label{relation-avec-Manin}
Every  homogeneous-quadratic homogeneous-symplectic Lie superalgebra $({\frak g}, B,\omega)$ over an algebraically closed field is a special homogeneous-symplectic homogeneous-Manin superalgebra $({\frak g}:={\frak a}\oplus {\frak b}, B,\omega)$. Moreover, we have either ${\frak z}({\frak g})\cap {\frak a}\neq \left\{0\right\}$ or ${\frak z}({\frak g})\cap {\frak b}\neq \left\{0\right\}$.\\
\end{thm}

\noindent\underline{\sl Proof} : In the first, we suppose that $({\frak g}, B,\omega)$ is an odd-quadratic odd-symplectic Lie superalgebra. In this case, the proof of the theorem is similar to the proof of Theorem 5.7 \cite{quadsymp}, where $({\frak g}, B,\omega)$ is a even-quadratic even-symplectic Lie superalgebra since the invariance of $B$ does not depend of the parity of $B$. Besides, the proof of Theorem \ref{relation-avec-Manin} when $({\frak g}, B,\omega)$ is an even-quadratic odd-symplectic Lie superalgebra is similar to the proof of Theorem \ref{relation-avec-Manin} when $({\frak g}, B,\omega)$ is odd-quadratic even-symplectic are similar, so we prove only Theorem \ref{relation-avec-Manin} when $({\frak g}, B,\omega)$ an even-quadratic odd-symplectic Lie superalgebra. Denotes $\Delta$ the unique odd invertible and skew-supersymmetric superderivation on ${\frak g}$ such that $\omega(X,Y)=B(\Delta(X),Y),\ \forall X,Y\in {\frak g}$ and Consider $\widetilde{\Delta}:= \left[\Delta,\Delta\right]$. Since $\widetilde{\Delta}$ is an even invertible and skew-symmetric superderivation on ${\frak g}$, then by the proof of Theorem 5.7 \cite{quadsymp}, we have $({\frak g}= {\frak a}\oplus {\frak b},B)$ is an even Manin superalgebra, where ${\frak a}:= \bigoplus_{\lambda \in Sp(\widetilde{\Delta})^+} {\frak g}^{\lambda}$ and ${\frak b}:= \bigoplus_{\lambda \in Sp(\widetilde{\Delta})^-} {\frak g}^{\lambda}$ (for more detail we can see \cite{quadsymp}). On the other hand, it is not difficult to show that ${\frak g}^{\lambda}$ is stable under $\Delta$, forall $\lambda$ in $Sp(\widetilde{\Delta})$. It comes that $\Delta({\frak a})\subseteq {\frak a}$ and $\Delta({\frak b})\subseteq {\frak b}$ and consequently, we obtain that $({\frak g}, B,\omega)$ is a special odd-symplectic even-Manin superalgebra. Similar to Theorem 5.7 \cite{quadsymp}, we have either ${\frak z}({\frak g})\cap {\frak a}\neq \left\{0\right\}$ or ${\frak z}({\frak g})\cap {\frak b}\neq \left\{0\right\}$.\ \ \ \ $\square$ \\

In the following, we are going to give inductive descriptions of special odd-symplectic odd-Manin superalgebras and inductive descriptions of odd (resp. even)-symplectic even (resp. odd)-Manin superalgebras.\\

\textbf{Inductive descriptions of special odd-symplectic odd-Manin superalgebras}
\\

Let $({\frak g}:= {\frak a}\oplus {\frak b}, B,\omega)$ be a special odd-symplectic odd-Manin superalgebra and $\delta$ the unique even invertible skew-supersymmetric superderivation on ${\frak g}$ such that $\omega(X,Y)= B(\delta(X),Y)$, ${\forall}_{X,Y\in {\frak g}}$ and $\delta({\frak a})\subseteq {\frak a}$ and $\delta({\frak b})\subseteq {\frak b}$. Suppose that $\bar{D}\in {(Der_s({\frak g}))}_{\bar{1}}$ which leaves ${\frak b}$ stable and $x_0\in {\frak b}_{\bar{0}}$ such that the conditions (\ref{condition-quad-odd})  are satisfied. Moreover, we suppose that there exist $\lambda \in {\mathbb K}\setminus \left\{0\right\}$ and $c_1\in {\frak b}_{\bar{1}}$ such that conditions (\ref{condition-oddquad-oddsymp}) are satisfied. Then the generalized odd double extension $(\widetilde{\frak g},\widetilde{B})$ defined as subsection 3.1 provided with bracket (\ref{bracket-quad-odd}) , with $k=0$ , the skew-supersymmetric bilinear form (\ref{form-quad-odd}) and the even invertible superderivation $\Delta$ defined as (\ref{derivation-quad-odd}) is a special odd-symplectic odd-Manin superalgebra which we call the generalized odd Manin double extension of $({\frak g}, B,\omega)$ be the one-dimensional Lie superalgebra with null even part by means of $(\bar{D}, x_0, c_1,\lambda)$.\\

\begin{prop}
Let $({\frak g}, B,\omega)$ be a special odd-symplectic odd-Manin superalgebra. Then ${\frak g}$ is a generalized odd double extension of a special odd-symplectic odd-Manin superalgebra $({\frak h}, B_{\frak h}, {\omega}_{\frak h})$ be the one-dimensional Lie superalgebra with null even part.\\
\end{prop}

\noindent\underline{\sl Proof} : Since ${\frak z}({\frak g})\cap {\frak g}_{\bar{0}}\neq \left\{0\right\}$, then  ${\frak z}({\frak g})\cap {\frak a}_{\bar{0}}\neq \left\{0\right\}$ or ${\frak z}({\frak g})\cap {\frak b}_{\bar{0}}\neq \left\{0\right\}$. We suppose that ${\frak z}({\frak g})\cap {\frak a}_{\bar{0}}\neq \left\{0\right\}$ and we denote $I:={\mathbb K}e^*$ where $e^*\in {\frak z}({\frak g})\cap {\frak a}_{\bar{0}}\setminus \left\{0\right\}$. Clearly, $I$ is totally isotropic and consequently $I\subseteq J$, where $J={I}^{\bot}$ with respect to $B$. From the non-degeneracy, the parity of $B$ and the fact that $B({\frak a},{\frak a})=\left\{0\right\}$, we deduce the existence of $e\in {\frak b}_{\bar{1}}$ such that $B(e^*,e)\neq 0$. Consider $W:= {\mathbb K}e^*\oplus {\mathbb K}e$. Since $W$ is non-degenerate, we have ${\frak g}= W\oplus {\frak h}$, where ${\frak h}$ is the orthogonal of $W$ with respect to $B$. It comes that $B_{\frak h}:= {B\vert}_{{\frak h}\times {\frak h}}$ is non degenerate and $J= I\oplus {\frak h}$. As ${\frak a}\subseteq J$, then ${\frak a}:= {\mathbb K}e^*\oplus \widetilde{\frak a}$, where $\widetilde{\frak a}= {\frak a}\cap {\frak h}$. Moreover, if we denote $\widetilde{\frak b}= {\frak b}\cap J= {\frak b}\cap {\frak h}$, then we obtain easily ${\frak h}= \widetilde{\frak a}\oplus \widetilde{\frak b}$. Now, by the proof of Proposition \ref{inverse-oddquad-oddsymp} , we have $({\frak h}, {\left[,\right]}_{\frak h}:= p\circ {\left[,\right]\vert}_{{\frak h}\times {\frak h}}, B_{\frak h}, {\omega}_{\frak h})$ is an odd-symplectic odd-quadratic Lie superalgebra, where $p: {\frak h}\oplus {\mathbb K}e^* \longrightarrow {\frak h}$ denote the projection $p(X+ \alpha e^*)= X$, for all $X\in {\frak h}$. The odd-symplectic form ${\omega}_{\frak h}$ on ${\frak h}$ is defined by ${\omega}_{\frak h}(X,Y)= B_{\frak h}(\delta(X),Y)$, $\forall X,Y\in {\frak h}$, where $\delta:= p\circ {\Delta\vert}_{\frak h}$. Once more following Proposition \ref{inverse-oddquad-oddsymp}, we have $({\frak g},B,\omega)$ is a generalized odd double extension of $({\frak h}, {\left[,\right]}_{\frak h}, B_{\frak h}, {\omega}_{\frak h})$ by means $(\bar{D},x_0,k)$ where $\bar{D}=p\circ {(\left[e,.\right])}_{\frak h}$ and $\left[e,e\right]= x_0 + ke^*$. Now it remains to verify that $\delta(\widetilde{\frak a})\subseteq \widetilde{\frak a}$, $\delta(\widetilde{\frak b})\subseteq \widetilde{\frak b}$, $\bar{D}(\widetilde{\frak b})\subseteq \widetilde{\frak b}$, $x_0\in {\widetilde{\frak b}}_{\bar{0}}$ and $k=0$. Since ${\frak b}$ is a sub-superalgebra of ${\frak g}$, then it is clear that $x_0\in {\widetilde{\frak a}}_{\bar{0}}$ and $k=0$. In addition, the fact that $J$ is a graded ideal of ${\frak g}$, ${\frak b}$ is a sub-superalgebra of ${\frak g}$ and $\Delta({\frak b})\subseteq {\frak b}$ and $\Delta({\frak a})\subseteq {\frak a}$ imply that $\delta$ and $\bar{D}$ leave stable $\widetilde{\frak b}$ and $\delta(\widetilde{\frak a})\subseteq \widetilde{\frak a}$ and so the result.\ \ \ \ $\square$ \\

\textbf{Inductive descriptions of special odd (resp. even)-symplectic even (resp. odd)-Manin superalgebras}
\\

In the following, we are going to introduce the notion of generalized double extension of a special odd (resp. even)-symplectic even (resp. odd)-Manin superalgebra by the two-dimensional abelian Lie supralgebra. This notion is a particular case of the notion of the generalized double extension of even (resp. odd)-quadratic odd (resp. even)-symplectic Lie superalgebras by the two-dimensional abelian Lie superalgebra.\\

 Let us consider $({\frak h}:= {\frak a}\oplus {\frak b}, B)$ an even (resp. odd)-Manin superalgebra, $D$ and  $\bar{D}$ two skew-supersymmetric superderivations of $({\frak h}, B)$ which are respectively even and odd and which both leave ${\frak b}$ stable, $x_0\in {\frak b}_{\bar{0}}$ and $x_1\in {\frak b}_{\bar{1}}$ such that the conditions (\ref{condition-even}) and (\ref{condition-compatibility-even}) of Theorem \ref{theorem-even} (resp. (\ref{condition-odd}) and (\ref{condition-compatibility}) of Theorem \ref{theorem}) are satisfied.
Then, the generalized (resp. odd) double extension $({\frak t}, \gamma)$ of $({\frak h},B)$ by the two-dimensional abelian Lie superalgebra by means of $(D,\bar{D},x_0,x_1,0)$, obtained in Theorem \ref{theorem-even} (resp. in Theorem \ref{theorem}) is an even (resp. odd)-Manin superalgebra. Explicitly ${\frak t}= \widetilde{\frak a}\oplus \widetilde{\frak b}$, where $\widetilde{\frak a}:= {\mathbb K}e_0^*\oplus {\frak a}\oplus {\mathbb K}e_1^*$ and $\widetilde{\frak b}:= {\mathbb K}e_0\oplus {\frak b}\oplus {\mathbb K}e_1$.\\

We call the even (resp. odd)-Manin superalgebra $({\frak t}= \widetilde{\frak a}\oplus \widetilde{\frak b},\gamma)$ the generalized (resp. odd) Manin double extension of $({\frak h}= {\frak a}\oplus {\frak b},B)$ by the two-dimensional Lie superalgebra ${\mathbb K}e_0\oplus {\mathbb K}e_1$ (by means of $(D,\bar{D}, x_0,x_1)$).\\

Now, we suppose that the even (resp. odd)-Manin superalgebra $({\frak h}= {\frak a}\oplus {\frak b}, B)$ admits a special odd (resp. even)-symplectic structure $\omega$. We denote  $\delta$ the unique odd invertible skew-supersymmetric superderivation on ${\frak h}$ such that $\omega(X,Y)= B(\delta(X),Y)$, $\delta({\frak a})\subseteq {\frak a}$ and $\delta({\frak b})\subseteq {\frak b}$. In addition, we consider $({\frak t}, \gamma)$ the generalized (resp. odd) Manin double extension of the even (resp. odd)-Manin superalgebra $({\frak h},B)$ by ${\mathbb K}e_0\oplus {\mathbb K}e_1$ by means of $(D, \bar{D}, x_0,x_1)$. If there exist $c_0\in {\frak b}_{\bar{0}}, \, c_1\in {\frak b}_{\bar{1}},$ such that the conditions (\ref{c1symp-impaire})-(\ref{c5symp-impaire}) of Proposition \ref{g-d-e-quad-paire-symp-impaire} (resp. (\ref{c1symp-paire})-(\ref{c5symp-paire}) of Proposition \ref{g-d-e-quad-impaire-symp-paire}) are satisfied, then the odd invertible skew-supersymmetric superderivation $\Delta$ defined as in Proposition \ref{g-d-e-quad-paire-symp-impaire} (resp. Proposition \ref{g-d-e-quad-impaire-symp-paire}), with $\alpha=0$, leaves stable $\widetilde{\frak b}$ and $\widetilde{\frak a}$. Consequently, $({\frak t}, \gamma, \Omega)$ is a special  odd (resp. even)-symplectic even (resp. odd)-Manin superalgebra such that $\Omega: {\frak t}\times {\frak t}\longrightarrow {\mathbb K}$ defined by
$$\Omega(X,Y):=\gamma(\Delta(X),Y), X,Y\in {\frak t}.$$

\begin{prop}
Let $({\frak g}= {\frak a}\oplus {\frak b},B,\omega)$ be a special odd (resp. even)-symplectic even (resp.odd)-Manin superalgebra and let $\Delta$ be the unique odd invertible skew-symmetric superderivation of ${\frak g}$ such that
$\omega(X,Y)= B(\Delta(X),Y)$, $\forall X,Y\in {\frak g}$, $\Delta({\frak a})\subseteq {\frak a}$ and $\Delta({\frak b})\subseteq {\frak b}$. If either ${\frak z}({\frak g})\cap {\frak a}\neq \left\{0\right\}$ or ${\frak z}({\frak g})\cap {\frak b}\neq \left\{0\right\}$, then ${\frak g}$ is a generalized (resp. odd) double extension of a special odd (resp. even)-symplectic even (resp. odd)-Manin superalgebra  $({\frak h}= \widetilde{\frak a}\oplus \widetilde{\frak b},B_{\frak h},{\omega}_{\frak h})$ by the two-dimensional abelian Lie superalgebra.\\
\end{prop}

\noindent\underline{\sl Proof} : We suppose that ${\frak z}({\frak g})\cap {\frak a}\neq \left\{0\right\}$ (the case where ${\frak z}({\frak g})\cap {\frak b}\neq \left\{0\right\}$ is similar). Since ${\mathbb K}$ is algebraically closed, $\Delta({\frak z}({\frak g})\cap {\frak a})\subseteq {\frak z}({\frak g})\cap {\frak a}$ and $\Delta$ is invertible, then there exist $e_0^*\in {({\frak z}({\frak g})\cap {\frak a})}_{\bar{0}}\setminus \left\{0\right\}$ and $e_1^*\in {({\frak z}({\frak g})\cap {\frak a})}_{\bar{1}}\setminus \left\{0\right\}$ such that $\Delta(e_0^*)=\lambda e_1^*$ and $\Delta(e_1^*)=\lambda e_0^*$ where $\lambda \in {\mathbb K}\setminus \left\{0\right\}$. We denote $I:= {\mathbb K}e_0^*\oplus {\mathbb K}e_1^*$ and $J$ its orthogonal with respect to $B$. As $B$ is even (resp. odd), non-degenerate and $B({\frak a}, {\frak a})=\left\{0\right\}$, then there exist $e_0\in {\frak b}_{\bar{0}}$ and  $e_1\in {\frak b}_{\bar{1}}$ such that $B(e_0^*,e_0)\neq \left\{0\right\}$ and $B(e_1^*,e_1)\neq \left\{0\right\}$ (resp. $B(e_0^*,e_1)\neq \left\{0\right\}$ and $B(e_1^*, e_0)\neq \left\{0\right\}$). Since the ${\mathbb Z}_2$-graded vector subspace $W:= I\oplus {\mathbb K}e_0\oplus {\mathbb K}e_1$ of ${\frak g}$ non-degenerate, we have ${\frak g}= W\oplus {\frak h}$, where ${\frak h}$ is the orthogonal of $W$ with respect to $B$ and it comes that $B_{\frak h}= {B\vert}_{{\frak h}\times {\frak h}}$ is non-degenerate. The fact that ${\frak a}\subseteq J$ and $J= I\oplus {\frak h}$, imply that ${\frak a}= I\oplus \widetilde{{\frak a}}$ where $\widetilde{{\frak a}}= {\frak a}\cap {\frak h}$. In addition, if we consider $\widetilde{{\frak b}}= {\frak b}\cap {\frak h}={\frak b}\cap J$, then we obtain ${\frak h}= \widetilde{{\frak a}}\oplus\widetilde{{\frak b}}$. Moreover $B_{\frak h}(\widetilde{\frak a},\widetilde{\frak a})= \left\{0\right\}= B_{\frak h}(\widetilde{\frak b}, \widetilde{\frak b})$.\\

\textbf{First case:} We suppose that $({\frak g},B,\omega)$ is a special odd-symplectic even-Manin superalgebra. Then, following Proposition \ref{inverse-even-quad-odd-symp-Lie} $({\frak h}, {\left[,\right]}_{\frak h}:={\left[,\right]\vert}_{{\frak h}\times {\frak h}}, B_{\frak h}, {\omega}_{\frak h})$ is an odd-symplectic even-quadratic Lie superalgebra. The odd-symplectic form ${\omega}_{\frak h}$ is defined by ${\omega}_{\frak h}(X,Y)= B_{\frak h}(\delta(X),Y), \forall X,Y\in {\frak h}$ such that $\delta = p_1\circ {\Delta\vert}_{\frak h}$ where $p_1: {\mathbb K}e_0^*\oplus {\mathbb K}e_1^* \oplus {\frak h}\longrightarrow {\frak h}$ denotes the projection $p_1(X+\alpha e_0^* + \beta e_1^*)= X, \forall X\in {\frak h}$. Once more by Proposition \ref{inverse-even-quad-odd-symp-Lie}, $({\frak g},B,\omega)$ is the generalized double extension of $({\frak h}, {\left[,\right]}_{\frak h}, B_{\frak h}, {\omega}_{\frak h})$ by ${\mathbb K}e_0\oplus {\mathbb K}e_1$ by means  $(D,\bar{D},x_0,x_1,c_0,c_1,\lambda)$, where $D= p_1\circ {ad_{\frak g}(e_0)\vert}_{\frak h}$, $\bar{D}=p_1\circ {ad_{\frak g}(e_1)\vert}_{\frak h}$, $ x_0 = p_2({\left[e_1,e_1\right]}_{\frak g})$, $ x_1 = p_2({\left[e_0,e_1\right]}_{\frak g})$, $c_0= p_2\circ \Delta(e_1)$ and $c_1= p_2\circ \Delta(e_0)$. where $p_2: {\frak g}\longrightarrow {\frak h}$ denotes the projection $p_2(\alpha e_0^*+ \beta e_1^*+X+\alpha'e_0 + \beta'e_1)= X$, for all  $X\in {\frak h}$, $\alpha, \beta, \alpha',\beta'\in {\mathbb K}$. Now, as $\widetilde{\frak a}= {\frak a}\cap {\frak h}$ and  $\Delta ({\frak a})\subseteq {\frak a}$, we deduce that $\delta (\widetilde{\frak a})\subseteq \widetilde{\frak a}$ (analogously, we have $\delta (\widetilde{\frak b})\subseteq \widetilde{\frak b}$). Consequently, we have $({\frak h}, {\left[,\right]\vert}_{{\frak h}\times {\frak h}}, B_{\frak h}, {\omega}_{\frak h})$ is a special odd-symplectic even-Manin superalgebra. On the other hand, since ${\frak b}$ is a sub-superalgebra of ${\frak g}$ and $e_0,e_1\in {\frak b}$, it follows that $x_0\in {\frak b}_{\bar{0}}$ and $x_1\in {\frak b}_{\bar{1}}$. Once more,  since ${\frak b}$ is a sub-superalgebra of ${\frak g}$ and $\widetilde{\frak b}= {\frak b}\cap {\frak h}$, it follows that $D$ and $\bar{D}$ leave stable $\widetilde{\frak b}$ and so the result.\\

\textbf{Second case:} We suppose that $({\frak g},B,\omega)$ is a special even-symplectic odd-Manin superalgebra. Using Proposition \ref{g-d-e-quad-impaire-symp-paire} and the same reasoning of the first case of this proof we obtain the result.\ \ \ \ $\square$ \\


\providecommand{\href}[2]{#2}

\end{document}